%% file: note.tex
\documentclass[final,onefignum,onetabnum]{siamart171218}

\usepackage{subfig}
\usepackage{xcolor} 
\usepackage{cancel} 
\usepackage{xspace} 

\newcommand{\vect}{\mathbf}

\newcommand{\ignore}[1]{}

\newcommand{\hypre}{{\em hypre}\xspace}
\graphicspath{{figures/}{/Users/sthomas1/Desktop/ecp_proposal_2016/milestone_reports/fy20_q2/figures/}}

\title{Two-Stage Gauss--Seidel Preconditioners and Smoothers for Krylov Solvers on a GPU cluster}

\author{
Stephen~Thomas%
\footnotemark[2]
\and Ichitaro~Yamazaki%
\footnotemark[1] 
\and Luc~Berger-Vergiat%
\thanks{Sandia National Laboratories, Albuquerque, New Mexico}
\and Brian~Kelley%
\footnotemark[1]
\and Jonathan~Hu%
\footnotemark[1]
\and Paul Mullowney%
\thanks{National Renewable Energy Laboratory, Golden, Colorado}
\and Sivasankaran~Rajamanickam%
\footnotemark[1]
\and Katarzyna~\'{S}wirydowicz%
\thanks{Pacific Northwest National Laboratory, Richland, WA}
}

\begin{document}

\maketitle

\begin{abstract}
Gauss-Seidel (GS) relaxation is often employed  as a preconditioner for a Krylov 
solver or as a smoother for Algebraic Multigrid (AMG).  
However, the requisite sparse triangular solve is difficult 
to parallelize on many-core architectures such as graphics processing units (GPUs).
In the present study, the performance of the sequential GS relaxation based on
a triangular solve is compared with two-stage variants,
replacing the direct triangular solve with a fixed number of inner Jacobi-Richardson (JR) iterations.
%
%
When a small number of inner iterations is sufficient to maintain the Krylov 
convergence rate, the two-stage GS (GS2) often outperforms the sequential algorithm
on many-core architectures.
The GS2 algorithm is also compared with JR.
When they perform the same number of flops for SpMV
(e.g. three JR sweeps compared to two GS sweeps with one inner JR sweep), 
the GS2 iterations, and the Krylov solver preconditioned with GS2,
may converge faster than the JR iterations.
Moreover, for some problems (e.g. elasticity), it was found that JR may diverge with a damping factor of one,
whereas two-stage GS may improve the convergence with more inner iterations.
Finally, to study the performance of the two-stage smoother and
preconditioner for a practical problem, 
these were applied to incompressible fluid flow simulations on GPUs.  
%
%
%
\end{abstract}

\section{Introduction}
\input{intro}

\section{Related Work}
\input{related}

\section{Gauss-Seidel Algorithms}
\label{sec:algo}

\input{algo}

\section{Two-stage Gauss-Seidel Algorithms}
\label{sec:algo-twostage}
\input{algo-twostage}

\section{Experimental Methodologies}
\label{sec:others}
\input{algo-others}

\section{Implementations}
\input{impl}

\section{Experimental Setup}
\input{result-setup}

\section{Experimental Results with Model Problems}
\label{sec:result-model}
\input{result-model}

\section{Experimental Results with Suite Sparse Matrices}

\input{result-ss}

\section{Experimental Results using Lower Precision Preconditioner}
\label{sec:result-mixed}
\input{result-mixed}

\section{Experimental Results with Nalu-Wind}
\label{sec:nalu-wind}
\input{result}

\section{Conclusion}
\input{conclusion}

\section*{Acknowledgment}
Funding was provided by the Exascale Computing Project (17-SC-20-SC).  The
National Renewable Energy Laboratory is operated by Alliance for Sustainable
Energy, LLC, for the U.S. Department of Energy (DOE) under Contract No.
DE-AC36-08GO28308.  Sandia National Laboratories is a multimission laboratory
managed and operated by National Technology \& Engineering Solutions of Sandia,
LLC, a wholly owned subsidiary of Honeywell International Inc., for the U.S.
Department of Energy National Nuclear Security Administration under contract
DE-NA0003525.
A portion of this research used resources of the Oak Ridge Leadership Computing
Facility, which is a DOE Office of Science User Facility supported under
Contract DE-AC05-00OR22725 and using computational resources at NREL,
sponsored by the DOE Office of Energy Efficiency and
Renewable Energy.
\bibliographystyle{elsarticle-num}
\bibliography{ref}
\end{document}

%% file: intro.tex
\label{sec:problem description}
Solving large sparse linear systems of the form $A\vect{x} = \vect{b}$ is a
basic and fundamental component of physics based modeling and simulation.  For
solving these large linear systems, Krylov methods such as the Generalized
Minimal Residual (GMRES)~\cite{Saad:1986} are widely adopted.  To
accelerate the convergence rate of the Krylov solver, Gauss-Seidel (GS)
relaxation is often employed either as the preconditioner for a Krylov solver
directly, or as the smoother in a $V$--cycle of a multigrid preconditioner.

In this study, GS relaxation is examined in the context of the
Trilinos~\cite{Heroux:2005} and \hypre~\cite{Falgout:2002} software
packages.  On a distributed-memory computer, both Trilinos and \hypre
implement a hybrid variant~\cite{Baker:2011} of Gauss-Seidel iteration where
the boundary information is exchanged to compute the residual vector, but then
each MPI process applies a fixed number of local relaxation sweeps
independently.  This hybrid block-Jacobi type algorithm
is shown to be effective and scalable for many problems (i.e., the iteration
counts stay roughly constant with increasing subdomain count).  
%
However, to perform the local relaxation, each process still requires a local
sparse-triangular solve, which is inherently challenging to
parallelize on the GPU architecture.  Although several techniques have been
proposed to improve the time to solution, the local sparse-triangular solve
recurrence could still be the bottleneck for scalability of the
solver, and in turn, of the entire simulation.  

In the present study, a two-stage iterative variant of the GS
relaxation algorithm is proposed as an alternative to the 
sequential algorithm based on a triangular solve. 
In the two-stage algorithm, henceforth called GS2, the triangular solve is replaced with 
a fixed number of  Jacobi-Richardson inner relaxation sweeps.  
Each inner sweep carries out  the same number of
floating-point operations (flops) as the triangular solve, but is based on
a sparse matrix vector multiply (SpMV), which is much easier to parallelize
and can be much faster than the sparse-triangular solve on a GPU.
Consequently, when a small number of inner iterations
is needed to maintain the quality of the preconditioner,
GS2 can surpass the sequential approach. 

The GS2 is also compared with the Jacobi-Richardson (JR) iteration,
which is a special case of the GS2 iteration without
inner iterations. Our results using Laplace and elasticity problems, 
together with the Suite-Sparse matrix collection,
demonstrate that when GS2 and JR are applied as stand-alone solvers
and execute the same number of SpMVs, then JR converges asymptotically faster.
Consequently,  GMRES often converges with fewer iterations, 
and with a shorter time to solution.
In contrast, the JR sweep computes the SpMV with the full matrix, while
the inner JR sweep of GS2 applies the SpMV with a triangular matrix. 
When the two methods execute about the same number of flops for SpMV
(e.g. three JR sweeps compared to two GS2 sweeps with one inner JR sweep),
the GS2 preconditioner may lead to faster convergence
as measured by the iteration counts, and in some cases, the time to
solution.  Moreover, for some problems (e.g. elasticity),
JR may diverge, 
while GS2 can improve the convergence using additional inner iterations
(if the inner iteration converges).

Lower precision floating point arithmetic is also investigated
as an alternative approximation for the preconditioner.
In experiments with a conjugate gradient (CG) solver and symmetric version of GS2 (SGS2),
high accuracy was not required for the outer SGS 
nor the inner JR iterations.  Thus, the low-precision inner 
iteration can be applied  to maintain the convergence rate of the 
Krylov solver, while reducing the compute time by a factor of $1.4\times$.

Finally, to study the effectiveness of the GS2 preconditioner
and smoother in a practical setting, they are applied in 
Nalu-Wind~\cite{Sprague:2020} for incompressible fluid flow simulations.  
In the context of a time-split,
pressure-projection scheme within nonlinear Picard iterations, 
GS is used as a preconditioner for GMRES when solving the
momentum problems, and GS, or Chebyshev, is applied as a smoother 
in algebraic multigrid (AMG) when solving pressure problems. For representative
wind-turbine incompressible flows, simulation results are reported for an
atmospheric boundary layer (ABL) precursor simulation on a structured mesh.
Our strong-scaling results 
on the Oak Ridge Leadership Computing Facility (OLCF) supercomputer Summit
demonstrate that the
GS2 preconditioner requires just one inner sweep to obtain a similar
convergence rate compared to the sequential algorithm, and thus it
achieves a faster time to solution. In addition, the
GS2 smoother was often as effective as a Chebyshev polynomial smoother.
%

%% file: related.tex
Because a sparse triangular solver is a critical kernel in many applications,
significant effort has been devoted to improving the performance of a
general-purpose sparse triangular solver for both 
CPUs~\cite{Anderson:1989,Saltz:1990,Bradley:2016} and
GPUs~\cite{Li:2013,Naumov:2011,Picciau:2016,Suchoski:2012,Li:2020}. 
For example, the traditional parallel algorithm is based
upon level-set scheduling, derived from the sparsity structure of the triangular
matrix and is also implemented in the cuSPARSE~\cite{cusparse} library.  
Independent computations proceed within each level of the tree.
In the context of  Gauss-Seidel relaxation, Deveci et.~al.
~implemented a multi-threaded triangular solver, based on multi-coloring \cite{Deveci:2016}.  
Multi-coloring exposes parallelism but may increase the iteration count for Krylov solvers in certain problems.
Another fine-grain algorithm is based on partial in-place inversion of products of elementary matrices after re-ordering to prevent fill-in. 
Here, the triangular solve is applied as a sequence of matrix-vector multiplies~\cite{Alvarado:1993}.
When the input matrices have a denser non-zero pattern arising from a 
super-nodal factorization, a block variant could exploit the GPU 
more effectively~\cite{yamazaki2020performance}.
%
Overall, the sparsity pattern of the triangular matrix may limit the amount of 
parallelism available for the solver to fully utilize many-core architectures,
especially when compared with other cuBLAS~\cite{cublas} and cuSPARSE operations, 
including a matrix-vector multiply.

Two-stage nested iteration, for which GS2 is a special case,
has been studied in a series of articles by Frommer and Szyld \cite{Szyld91,Szyld92a,Szyld92b,Szyld94}. 
Unlike the triangular solver algorithm, the two-stage iteration
can be implemented using only (sparse or dense) matrix-vector
products and vector updates. As stated previously, 
these operations can be implemented very
efficiently on multi-core architectures.
The application of Jacobi iterations to 
solve sparse triangular linear systems for ILU
preconditioners, rather than forward or backward recurrence, 
was also recently proposed by Chow et. al.~\cite{ChAnScDo2018}.

%% file: algo.tex
To solve a linear system $A \vect{x} = \vect{b}$, 
the Gauss-Seidel (GS) iteration is
based on the matrix splitting $A = L + D + U$, where $L$ and $U$ are the strictly
lower and upper triangular parts of the matrix~$A$, respectively.
Then, the sequential GS iteration updates the solution based on the following recurrence,
\begin{equation}\label{eq:one-stage}
 \vect{x}_{k+1} := \vect{x}_k + M^{-1} \vect{r}_k,  \quad k=0, 1, 2, \ldots
\end{equation}
where $\vect{r}_k = \vect{b} - A \vect{x}_k$, and $M = L + D$, $N=-U$ or $M = U + D$ 
for the forward or backward sweeps, respectively. The iteration matrix  
is $G = I-M^{-1}A$ and when (\ref{eq:one-stage}) converges,
$\rho(G) < 1$, (c.f. Theorem 4.1 of Saad~\cite{Saad03}).
When the matrix $A$ is symmetric, a symmetric variant of GS (SGS)
performs the forward followed by the backward sweep
to maintain the symmetry of the matrix operation.
In the remaining discussions, we refer to \ref{eq:one-stage} as the {\em non-compact form}.

To avoid explicitly forming the inverse $M^{-1}$ in 
\eqref{eq:one-stage}, a sparse triangular solve is applied
to the current residual vector~$\vect{r}_k$. On a distributed-memory 
computer, both Trilinos and \hypre 
distribute the matrix and the vectors in a 1--D
block row fashion among the MPI processes 
(the local matrix on the $p$-th MPI process is
rectangular, consisting of the square diagonal block~$A^{(p)}$ 
for the rows owned by the
MPI process and the off-diagonal block $E^{(p)}$ for the off-process columns.
On each MPI rank, \hypre stores the local diagonal block $A^{(p)}$
separately from the off-diagonal block $E^{(p)}$, whereas Trilinos stores them as a single matrix).
The standard parallel sparse triangular solve is based
on level-set scheduling~\cite{Anderson:1989,Saltz:1990}
and will compute the independent set of the solution elements
in parallel at each level of scheduling. 
Unfortunately, the sparsity structure of the triangular matrix may limit the parallelism that
the solver can exploit (e.g. the sparsity structure may lead to a long chain
of dependencies with a small number of solution elements that can be computed at
each level).  In addition, at the start of each level, neighboring processes
need to exchange the elements of the solution vector on a 
processor boundary for updating their local right-hand-side vectors.

To improve scalability, both Trilinos and \hypre implement a
hybrid variant of GS relaxation~\cite{Baker:2011}, where the neighboring 
processes first exchange the elements of the solution vector on the boundary, 
but then each process independently applies the local relaxation.
Hence, $M$ is a block diagonal matrix with each diagonal block $M^{(p)}$ corresponding to the
triangular part of the local diagonal block~$A^{(p)}$ on each process. 
The method may be considered as an inexact block Jacobi iteration, 
where the local problem associated with the diagonal block is solved using a single GS iteration.
Furthermore, in \hypre, each process may apply the multiple local GS sweeps
for each round of the neighborhood communication. With this approach,
each local relaxation updates only the local part of the vector $\vect{x}_{k+1}$
(during the local relaxation, the non-local solution elements on the boundary 
are not kept consistent among the neighboring processes).
Hence, the local GS relaxation, on the $p$-th MPI rank,
iterates to solve the local linear system $A^{(p)} \vect{x}^{(p)}_{k+1} = \widehat{\vect{b}}_t^{(p)}$
where 
$\widehat{\vect{b}}_t^{(p)} := \vect{b}^{(p)} - E^{(p)} \vect{y}_t^{(p)}$ and
$\vect{y}_t^{(p)}$ is the non-local part of the current solution vector exchanged
before the start of the local GS iteration.
In contrast, Trilinos currently computes the global residual vector for each local 
iteration through neighborhood communication.  The only exception is the symmetric GS 
relaxation that applies both the forward and backward
iterations after a single round of neighborhood communication. 

In this hybrid GS, each process can apply the local relaxation independently, 
improving the scalability compared to the global sparse triangular 
solve. This hybrid algorithm is shown
to be effective and scalable for many problems~\cite{Baker:2011} 
(e.g. the Krylov iteration count, when combined with this hybrid GS preconditioner, 
remains roughly constant with an increasing process count).
However, to implement the local iteration, each process  must still perform a 
local sparse triangular solve, which can be
difficult to parallelize on many-core architectures.

%% file: algo-twostage.tex
In the present study, a two-stage GS relaxation, with a fixed number of
``inner'' stationary iterations, approximately solves the triangular system with~$M$,
\begin{equation}\label{eq:two-stage}
 \widehat{\vect{x}}_{k+1} := \widehat{\vect{x}}_k + \widehat{M}^{-1} (
\vect{b} - A \widehat{\vect{x}}_k ), \quad k=0, 1, 2, \ldots
\end{equation}
where $\widehat{M}^{-1}$ represents the approximate triangular system solution,
i.e. $\widehat{M}^{-1} \approx M^{-1}$. 

A Jacobi-Richardson (JR) iteration is employed for the inner iteration.
More specifically, $\vect{g}_k^{(j)}$ denotes the approximate
solution from the $j$--th inner JR iteration at the $k$--th outer GS iteration. 
The initial solution is chosen to be the diagonally-scaled residual vector,
\begin{equation}\label{eq:jr-initial-guess}
  \vect{g}_k^{(0)} = D^{-1}\vect{r}_k,
\end{equation}
and the $(j+1)$--st JR iteration computes the solution by the recurrence
\begin{eqnarray}
\label{eq:jr}
  \vect{g}_k^{(j+1)} & := & \vect{g}_k^{(j)} + D^{-1} (\vect{r}_k - (L+D) \vect{g}_k^{(j)})\\
\label{eq:jacobi}
                     &  = & D^{-1}(\vect{r}_k - L \vect{g}_k^{(j)}).
\end{eqnarray}
Figure~\ref{code:two-stage} displays the SGS2 algorithm, the symmetric extension of GS2.  

\begin{figure}
\centerline{
 \fbox{\begin{minipage}[t]{.5\textwidth}
 {\normalsize\texttt{\input{./codes/gs2}}}
 \end{minipage}}
}
\caption{\label{code:two-stage}
Pseudo-code of two-stage hybrid Symmetric Gauss-Seidel iteration, using Jacobi-Richardson as the inner sweep.  $n_t$ is the number of applications of GS2.  $n_k$ is the number of local outer sweeps, and $n_j$ is the number of local inner sweeps.}
\end{figure}

When ``zero'' inner iterations are performed, the GS2 recurrence becomes
\[
  \widehat{\vect{x}}_{k+1} := \widehat{\vect{x}}_k + \vect{g}_k^{(0)} = 
  \widehat{\vect{x}}_k + D^{-1}(\vect{b} - A \widehat{\vect{x}}_k),
\]
and this special case corresponds to Jacobi-Richardson iteration on the global system, or local system on each process. When $s$ inner iterations are performed, it is easy to see
that
\begin{align*}
\widehat{\vect{x}}_{k+1} &:= \widehat{\vect{x}}_k + \vect{g}_k^{(s)} = \widehat{\vect{x}}_k + \sum_{j=0}^{s}
(-D^{-1} L)^j D^{-1} \widehat{\vect{r}}_k \\
&\approx \widehat{\vect{x}}_k + (I+D^{-1} L)^{-1} D^{-1} \widehat{\vect{r}}_k = \widehat{\vect{x}}_k + M^{-1}  \widehat{\vect{r}}_k,
\end{align*}
where $M^{-1}$ is approximated by the degree-$s$
Neumann expansion.
Note that $D^{-1} L$ is strictly lower triangular and thus the Neumann series 
converge in a finite number of steps.

\paragraph{Compact form of Gauss-Seidel recurrence}
The GS2 recurrence~\eqref{eq:two-stage} may be written as
\begin{eqnarray}
\label{eq:two-stage-1}
 \widehat{\vect{x}}_{k+1} & := & \widehat{\vect{x}}_k + \
 \widehat{M}^{-1} ( \vect{b} - (M - N) \widehat{\vect{x}}_k )\\
\label{eq:two-stage-2}
			  &  = & (I - \widehat{M}^{-1}M) \widehat{\vect{x}}_k
+ \widehat{M}^{-1} ( \vect{b} + N \widehat{\vect{x}}_k).
\end{eqnarray}
In the classical one-stage recurrence~\eqref{eq:one-stage}, the 
preconditioner matrix is taken as
$\widehat{M}^{-1} = M^{-1}$, and only the second term remains in the
recurrence~\eqref{eq:two-stage-2}, leading to the following ``compact'' form,
\begin{eqnarray}\label{eq:one-stage-2} 
\vect{x}_{k+1} &  := & M^{-1} (
\vect{b} + N \vect{x}_k).  
\end{eqnarray}
Hence, the recurrences~\eqref{eq:one-stage} and \eqref{eq:one-stage-2} are
mathematically equivalent, while the recurrence~\eqref{eq:one-stage-2} has a
lower computation cost.
The multi-threaded GS implementation in Trilinos is based on the compact recurrence.

A similar ``compact'' recurrence for the GS2 algorithm may be derived as
\begin{eqnarray}\label{eq:two-stage-3}
 \widetilde{\vect{x}}_{k+1} & := & \widehat{M}^{-1} ( \vect{b} + N \widetilde{\vect{x}}_k).
\end{eqnarray}
However, with the approximate solution using $\widehat{M}^{-1}$, the
recurrences~\eqref{eq:two-stage} and \eqref{eq:two-stage-3} are no longer
equivalent.  For example, even if it is assumed that $\widehat{\vect{x}}_k = \vect{x}_k$, comparing \eqref{eq:one-stage} and \eqref{eq:two-stage},
the difference in the residual norms using the classical and the standard
two-stage iterations is given by
\begin{eqnarray}
\nonumber
 \|\widetilde{\vect{r}}_{k+1} - \vect{r}_{k+1}\| & =   &\|A(I - \widehat{M}^{-1}M) M^{-1}\vect{r}_k\|\\
\label{eq:rnorm1}
& \le & \|A (I - \widehat{M}^{-1}M)\| \|M^{-1}\vect{r}_k\|
\end{eqnarray}
On the other hand, even assuming that $\widetilde{\vect{x}}_k = \vect{x}_k$,
comparing \eqref{eq:one-stage-2} and \eqref{eq:two-stage-3},
the difference between the classical and the compact two-stage relaxation is
\begin{eqnarray}
\nonumber
 \|\widehat{\vect{r}}_{k+1} - \vect{r}_{k+1}\| 
  & =   &\|A(I - \widehat{M}^{-1}M) (M^{-1}\vect{r}_k + \vect{x}_k)\|\\
\label{eq:rnorm2}
  & \le & \|A (I - \widehat{M}^{-1}M)\| \|M^{-1}\vect{r}_k\| + \|A (I - \widehat{M}^{-1}M)\| \|\vect{x}_k\|
\end{eqnarray}
and the compact form has the extra term with $\|\vect{x}_k\|$ in the bound.
%
%
For the recurrence~\eqref{eq:two-stage-3} to be as effective as the
recurrence~\eqref{eq:two-stage-1}, it was found that 
additional inner iterations are often required
(to make $\|I - \widehat{M}^{-1}M\|$ small).

\paragraph{Damping factor}

The convergence rate of the Gauss-Seidel iteration may be improved using
a matrix splitting parameterized with a damping factor $\omega$:
\begin{eqnarray}
\nonumber
\vect{x}_{k+1} & := & \vect{x}_k + \omega (D + \omega L)^{-1} (\vect{b} - A \vect{x}_k)\\
\nonumber
               &  = & \omega(D + \omega L)^{-1} (\vect{b} - [U + (1 - \frac{1}{\omega}) D] \vect{x}_k)
\end{eqnarray}
for the non-compact and compact forms, respectively.  
Therefore, the inner Jacobi iteration becomes
\begin{eqnarray}
\nonumber
 \vect{g}_k^{(j+1)} & := & \vect{g}_k^{(j)} + D^{-1}(\vect{r}_k - [\omega L + D] \vect{g}_k^{(j)})\\
\nonumber
                    &  = & D^{-1}(\vect{r}_k - \omega L \vect{g}_k^{(j)})
\end{eqnarray}
Moreover, for GS2, the convergence of the inner 
JR iteration may be improved by using another damping factor $\gamma$:
\begin{eqnarray}
\nonumber
 \vect{g}_k^{(j+1)} & := & \vect{g}_k^{(j)} + \gamma D^{-1}(\vect{r}_k - [\omega L + D] \vect{g}_k^{(j)})\\
\nonumber
                &  = & \gamma D^{-1}(\vect{r}_k - [\omega L + (1 - \frac{1}{\gamma}) D] \vect{g}_k^{(j)})\\
\nonumber
                &  = & (1 - \gamma) \vect{g}_k + \gamma D^{-1}(\vect{r}_k - \omega L \vect{g}_k)
\end{eqnarray}


%% file: codes/gs2.tex
\begin{tabbing}
m\=mmm\=mm\=mm\=m\=m\=m\=\kill
\> {\bf for} $t=1,2,\dots, n_t$ {\bf do}\\
\> 1. \> \emph{// exchange interface elements of current solution}\\
\>\> {\bf for} $k=1,2,\dots, n_k$ {\bf do}\\
\> 2. \>\> \emph{// compute new residual vector for forward sweep}\\
\>    \>\> $\vect{r}_k := \vect{b} - U \vect{x}_k$\\
\> 3. \>\> \emph{// perform local inner Jacobi iteration}\\
\>    \>\> $\vect{g}^{(0)}_k := D^{-1}\vect{r}_k$\\
\>    \>\> {\bf for} $j=0,1,\dots, n_j-1$ {\bf do}\\
\>    \>\>\> $\vect{g}^{(j+1)}_{k} := D^{-1}(\vect{r}_k - L \vect{g}^{(j)}_k)$\\
\>    \>\> {\bf end for}\\
\> 4. \>\> \emph{// update solution vector}\\
\>    \>\> $\vect{x}_{k+1} := \vect{x}_k + \vect{g}^{(s)}_{k}$\\
\> 5. \>\> \emph{// compute new residual vector for backward sweep}\\
\>    \>\> $\vect{r}_k := \vect{b} - L \vect{x}_k$\\
\> 6. \>\> \emph{// perform local inner Jacobi iteration}\\
\>    \>\> $\vect{g}^{(0)}_k := D^{-1}\vect{r}_k$\\
\>    \>\> {\bf for} $j=0,1,\dots, s$ {\bf do}\\
\>    \>\>\> $\vect{g}^{(j+1)}_{k} := D^{-1}(\vect{r}_k - U \vect{g}^{(j)}_k)$\\
\>    \>\> {\bf end for}\\
\> 7. \>\> \emph{// update solution vector}\\
\>    \>\> $\vect{x}_{k+1} := \vect{x}_k + \vect{g}^{(s)}_{k}$\\
\>    \> {\bf end for}\\
\> {\bf end for}
\end{tabbing}

%% file: algo-others.tex

Here, different variants of GS2 are compared, including JR, and 
numerical experiments are presented in Section~\ref{sec:result-model}.

{\it Comparison of GS2 with JR
when they perform the same number of sparse-matrix vector kernels (latency)}:
One sweep of GS2 with one inner JR sweep performs two SpMVs
(one SpMV with $A$ and another with~$L$).
\begin{eqnarray}
\label{eq:rk_gs2}
\vect{r}^{\mbox{gs2}(1)}_{k+1} & := & (I - AD^{-1}[I-LD^{-1}])^{k+1} \vect{r}_0\\
\nonumber
                               &  = & (I - AD^{-1} + AD^{-1}LD^{-1})^{k+1} \vect{r}_0
\end{eqnarray}
For a comparison, two sweeps of JR, which calls SpMV twice with $A$,
is given by the following recurrence:
\begin{equation}
\vect{x}_{k+2}^{\mbox{jr}} = D^{-1}\vect{r}^{\mbox{jr}}_k - D^{-1} A D^{-1}\vect{r}^{\mbox{jr}}_k,
\end{equation}
and
\begin{eqnarray}
\label{eq:rk_jr}
\vect{r}_{2k}^{\mbox{jr}} & := & (I - AD^{-1})^{2k} \vect{r}_0\\
\nonumber
                          &  = & (I - AD^{-1}[I - (L + U)D^{-1}])^{k} \vect{r}_0\\
\nonumber
                          &  = & (I - AD^{-1} + AD^{-1} (L + U)D^{-1})^{k} \vect{r}_0
\end{eqnarray}

In~\cite{Szyld92a}, it is shown that
when the matrix splitting $A = M-N$ is regular (i.e. $M^{-1} \ge 0$ and $N \ge 0$)
and $M = B - C$ is weakly regular (i.e. $M^{-1} \ge 0$ and $M^{-1}N \ge 0$), 
then $\rho(T_{n_j+1}) \ge \rho(T_1)^{n_j+1}$
where $T_{n_j+1}$ is the iteration matrix of the two-stage method with $n_j$ inner iterations. GS2 is a special case of this result.
$T_{1}$ corresponds to the JR iteration.
In~\cite{Szyld92a}, $\vect{g}_0 := D^{-1}\vect{r}_k$ is taken as the approximate solution after 
one inner sweep, whereas the initial guess~\eqref{eq:jr-initial-guess} 
is used here.
Furthermore, with no inner sweep (i.e. $n_j=0$), two sweeps of JR converges asymptotically faster than
one sweep of GS2 (i.e. $\rho(T_{2}) \ge \rho(T_1)^{2}$).

{\it Comparison of GS2 with JR
when they perform the same number of flops for 
the sparse-matrix vector kernels (data pass)}:
Because one inner sweep of GS2 applies an SpMV with the triangular matrix,
it requires about half the number of flops compared to the SpMV for one sweep of JR.
Hence, three sweeps of JR perform about the same number of flops for the SpMV
as the two sweeps of GS2 or one sweep of SGS2:

The SGS2 recurrence is given by
\begin{eqnarray}
\nonumber
\vect{r}^{\mbox{sgs2}(1)}_{k+1} & := & (I - AD^{-1}[I-UD^{-1}]) (I - AD^{-1}[I-LD^{-1}]) \vect{r}_k
\end{eqnarray}
while three sweeps of JR results in
\begin{eqnarray}
\nonumber
\vect{r}_{k+3}^{\mbox{jr}} & := & (I - AD^{-1})^{3} \vect{r}_k
\end{eqnarray}
According to the results in~\cite{Szyld92a},
two sweeps of GS2 converges asymptotically slower than four sweeps
of JR (i.e. $\rho(T_{2})^2 \ge \rho(T_1)^{4}$).
However, even for this special case of a regular splitting 
or  weak regular splitting, it may occur that $\rho(T_{2})^2 < \rho(T_1)^{3}$.
In such cases, GMRES with two-sweeps of the GS2 preconditioner
may converge faster than with three sweeps of the JR preconditioner.



{\it Comparison of GS2 with one inner sweep with GS2 in compact form with two inner sweeps}:
In Section~\ref{sec:algo} (before the discussion on the damping factor), it was
noted that the compact form often requires additional inner sweeps to match the
quality of the preconditioner using the non-compact form. Thus, 
the compact form using two inner sweeps is compared with
the non-compact using one inner sweep because these perform
the same number of flops for an SpMV.



%% file: impl.tex
Our experiments were conducted using Trilinos or \hypre.
In this section,
the implementations are described, which are now available in 
these two software packages.

\subsection{Trilinos}
\input{impl-trilinos}

\subsection{\hypre}

\input{impl-hypre}

%% file: impl-trilinos.tex
Kokkos Kernels provides performance portable
sparse and dense linear algebra and graph algorithms.
It now includes several implementations of the Gauss-Seidel
iteration, all of which are available for Belos solvers through the Ifpack2
interface (Belos and Ifpack2 are the Trilinos software packages
that provide iterative linear solvers and
algebraic preconditioners/smoothers, respectively). 
Brief descriptions of the implementations are provided here.

\subsubsection{Sequential Gauss-Seidel} \label{sequentialGS}

Based on the compact form of the recurrence~\eqref{eq:one-stage-2}, the
sequential implementation of Gauss-Seidel in Trilinos computes the solution 
from the first element to the last based on the natural ordering of the matrix.

\medskip
\fbox{\begin{minipage}[t]{.7\textwidth}
\begin{tabbing}
m\=m\=m\=m\=m\=m\=m\=\kill
\> {\bf for} $i=1,2,\dots, n$ {\bf do}\\
\>\> $\mathit{sum} := 0$ \\
\>\> {\bf for} each nonzero off-diagonal entry $a_{i,j}$ in the $i$th row of $A$ {\bf do}\\
\>\>\> $\mathit{sum} := \mathit{sum} + a_{i,j} x_{j}$\\
\>\> {\bf end for}\\
\>\> $x_i := (b_{i} - \mathit{sum})/a_{i,i}$\\
\> {\bf end for}
\end{tabbing}
\end{minipage}}
\medskip

The above algorithm is employed as the baseline implementation in comparisons,
in particular for numerical studies (e.g. convergence rates).
When the natural ordering of the matrix respects the flow direction of the
underlying physics, then this implementation results in the ``optimal''
convergence with the Gauss-Seidel preconditioner \cite{Elman14115}.

\subsubsection{Multi-threaded Gauss-Seidel}

This is an implementation of multicolor Gauss-Seidel in Kokkos Kernels \cite{Deveci:2016}. 
It is abbreviated as MT in this paper. It computes a parallel greedy coloring of the matrix, 
so that if two rows $i,j$ have the same color, then $a_{i,j}=a_{j,i}=0$. This means that 
the procedure in Section~\ref{sequentialGS} is parallel, without data races or 
nondeterministic behavior.

\medskip
\fbox{\begin{minipage}[t]{.7\textwidth}
\begin{tabbing}
m\=m\=m\=m\=m\=m\=m\=\kill
\> {\bf for} each color $c$ {\bf do}\\
\>\> {\bf for} each row $i$ with color $c$ {\bf do in parallel}\\
\>\>\> $\mathit{sum} := 0$ \\
\>\>\> {\bf for} each nonzero off-diagonal entry $a_{i,j}$ in the $i$th row of $A$ {\bf do}\\
\>\>\>\> $\triangleright$ $x_{j}$ is not modified during this parallel loop\\
\>\>\>\> $sum := sum + a_{i,j} x_{j}$\\
\>\>\> {\bf end for}\\
\>\>\> $x_i := (b_i - \mathit{sum})/a_{i,i}$\\
\>\> {\bf end for}\\
\> {\bf end for}
\end{tabbing}
\end{minipage}}
\medskip

The backward sweep is identical except it iterates over the colors in reverse. 
The number of colors required is bounded above by the maximum degree of the 
matrix~\cite{Deveci:2016}. Unlike sequential GS, GS-MT does not respect indirect 
dependencies between rows, so GS-MT usually increases the iteration count.
However, parallel execution still results in a much lower time to solution 
for sufficiently large matrices (e.g. see Fig.~\ref{tab:cg-default}).

\subsubsection{Two-stage Gauss-Seidel}

Instead of implementing a specialized Kokkos kernel,
the current implementation of GS2
relies on Kokkos Kernels for the required matrix or vector operations  
(e.g. SpMV).
The resulting code design is modular and thus allows access to the SpMV
kernels whose vendor-optimized implementations are often available for the
specific node architecture (e.g. MKL on Intel CPUs or cuSPARSE on NVIDIA GPUs).
However, this launches multiple kernels
and requires extra memory to explicitly store the lower or upper
triangular matrices, $L$ and $U$.  Because the lower-triangular matrices are
explicitly stored, separate from the original matrix $A$, the triangular matrix
is prescaled with the diagonal matrix ($D^{-1}L$ is stored) to avoid the extra
kernel launch for scaling the vector at each inner sweep.
Because achieving high performance for the two-stage iteration relies upon the
SpMV kernel, $A$ may be reordered.  Because the natural ordering often
provides rapid convergence, the natural ordering of $A$ is specified
for our current implementation.

%% file: impl-hypre.tex
\hypre provides several different Gauss-Seidel smoother algorithms. These
include symmetric GS, hybrid $\ell_1$ GS forward  and hybrid $\ell_1$ 
GS backward solve which perform a single round of neighborhood
communication followed by local solve(s)~\cite{Baker:2011}. 
The default smoother is hybrid symmetric GS. Initially, an
$\ell_1$ Jacobi smoother was available for the GPU. 
A GS2 iteration was implemented as a part of the 
\hypre solver stack as well. \hypre does not use Kokkos and
the GPU-accelerated components are implemented directly in CUDA. 
Each CUDA software release comes
    with highly optimized cuBLAS and cuSPARSE libraries that provide all basic
linear algebra operations (such as SpMV), however, custom kernels were
implemented to minimize the total number of kernel launches required.

\ignore{
For example, the JR implementation in \hypre is based on:
\begin{eqnarray}
\vect{x}_{k+1} & = & \vect{x}_k + (D^{(p)})^{-1}\vect{r}_{k} +  \omega\: (D^{(p)})^{-1}\left[\,
                     \vect{r}_{k} - A^{(p)} \: (D^{(p)})^{-1}\vect{r}_{k} \,\right]\\
               & = & \vect{x}_k + \omega\: (D^{(p)})^{-1}\left[\,
                     \vect{r}_{k} - (A^{(p)} - \frac{1}{\omega}D^{(p)}) \: (D^{(p)})^{-1}\vect{r}_{k} \,\right].
\end{eqnarray}
Moreover, after the two JR sweeps with $\vect{x}_0 = \vect{0}$, we have
\begin{eqnarray}
\vect{x}_{2} & = & (D^{(p)})^{-1}\vect{b} +  \omega\: (D^{(p)})^{-1}\left[\,
                     \vect{b} - A^{(p)} \: (D^{(p)})^{-1}\vect{b} \,\right]\\
             & = & \omega\: (D^{(p)})^{-1}\left[\,
                   \vect{b} - (A^{(p)} - \frac{1}{\omega}D^{(p)}) \: (D^{(p)})^{-1}\vect{b} \,\right].
\end{eqnarray}
}

All parallel matrices are stored in CSR-like format with two parts: 
the {\it diagonal part}  (which is local to the rank) and the {\it off-diagonal} part
that is used to multiply the data after exchange between ranks. 
Thus, the \hypre two-stage iteration becomes
\begin{equation}
\widehat{\vect{x}}_{k+1} :=  \widehat{\vect{x}}_k + \sum_{j=0}^{s}
(-D^{-1} L)^j D^{-1} \widehat{\vect{r}}_k
\label{eq:hypresmoother}
\end{equation}
%
%
%
where $ A^{\text{diag}} = D + L + U$ is the diagonal (on rank) part of $A$. In
the general case, the matrix $ A^{\text{diag}}$ is not triangular.

Computing $\vect{r}_k=\vect{b}-A\vect{x}_k$ in \hypre  requires a global SpMV: this means, two
SpMVs are needed, one with the local square matrix $A^{\text{diag}} $ and
another with the off-rank columns $A^{\text{offdiag}} $ and then the resulting
vectors must be added together.  Thus, the operations implied by
~(\ref{eq:hypresmoother}) are implemented as two CUDA kernels. In the first
kernel, the residual $\vect{r}_k$ is computed using two SpMVs. The kernel returns both
$\vect{r}_k$ and $D^{-1}\vect{r}_k$. The second kernel takes care of the remaining operations.

This particular approach requires additional storage for temporary results,
such as, $D^{-1}\vect{r}_k$. In order to optimize the implementation, temporary
vectors are allocated only once at each level and persist as long as $A$ remains
in the GPU global memory, which removes a large quantity of costly \verb|cudaMalloc| and \verb|cudaFree| calls.
Data exchange between ranks is required in our approach; but after the data is
exchanged once, all other operations are local to the rank.

%% file: result-setup.tex
Experiments were conducted on the Summit supercomputer at Oak Ridge National
Laboratory. Each node of Summit has two 22-core IBM Power9 CPUs and six NVIDIA
Volta 100 GPUs.  
For the solver performance studies in Sections~\ref{sec:result-model} to \ref{sec:result-mixed}, the development branch of Trilinos was employed
and for Section~\ref{sec:nalu-wind}, 
the master branch of Nalu-Wind was used. 
%
The code is compiled using CUDA version 10.1.168 and GNU {\tt gcc} version 7.4.0.

Except for the experiments with Nalu-Wind,
a random vector $\vect{b}$ is generated for the linear system,
and the iterative solver has converged when the relative 
residual norm is reduced by nine orders of magnitude.
In addition, unless explicitly specified,
a single GPU is employed to study convergence and performance.
For our discussion, GS($n_t$, $n_k$) and GS2($n_t$, $n_k$, $n_j$) refer
to the classical and two-stage GS, respectively, while JR($n_t$, $n_k$)
is Jacobi-Richardson,
where $n_t$ and $n_k$ are the global and local GS or JR sweeps, 
respectively, while $n_j$ is the number of inner JR sweeps
When the default $n_k = 1$ is employed, these are referenced
without $n_k$, e.g. GS2($n_t$, $n_j$).
Unless explicitly specified, the non-compact form of GS2
is applied.

%% file: result-model.tex
The convergence of GS2 and JR iterations are now compared for
Laplace and elasticity problems using Trilinos 
(these model problems are available through Xpetra Galeri).

\begin{figure}
\centerline{
  \includegraphics[width=.5\textwidth]{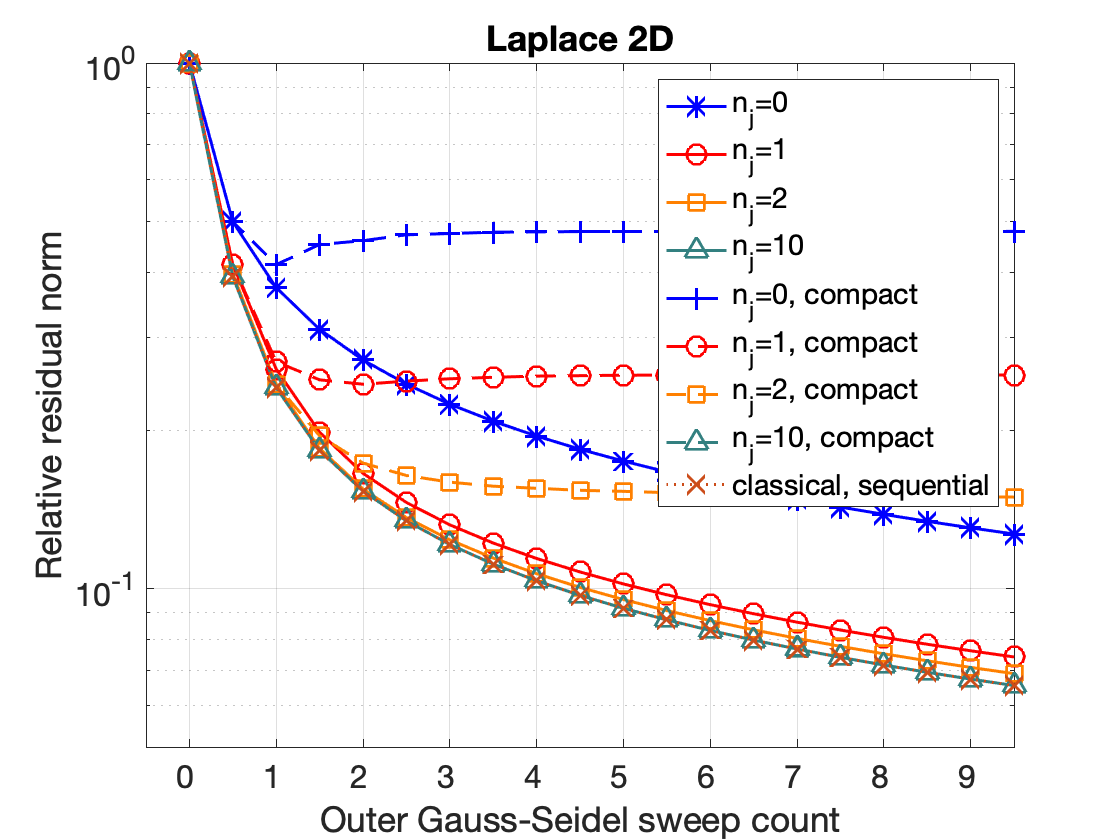}
  \includegraphics[width=.5\textwidth]{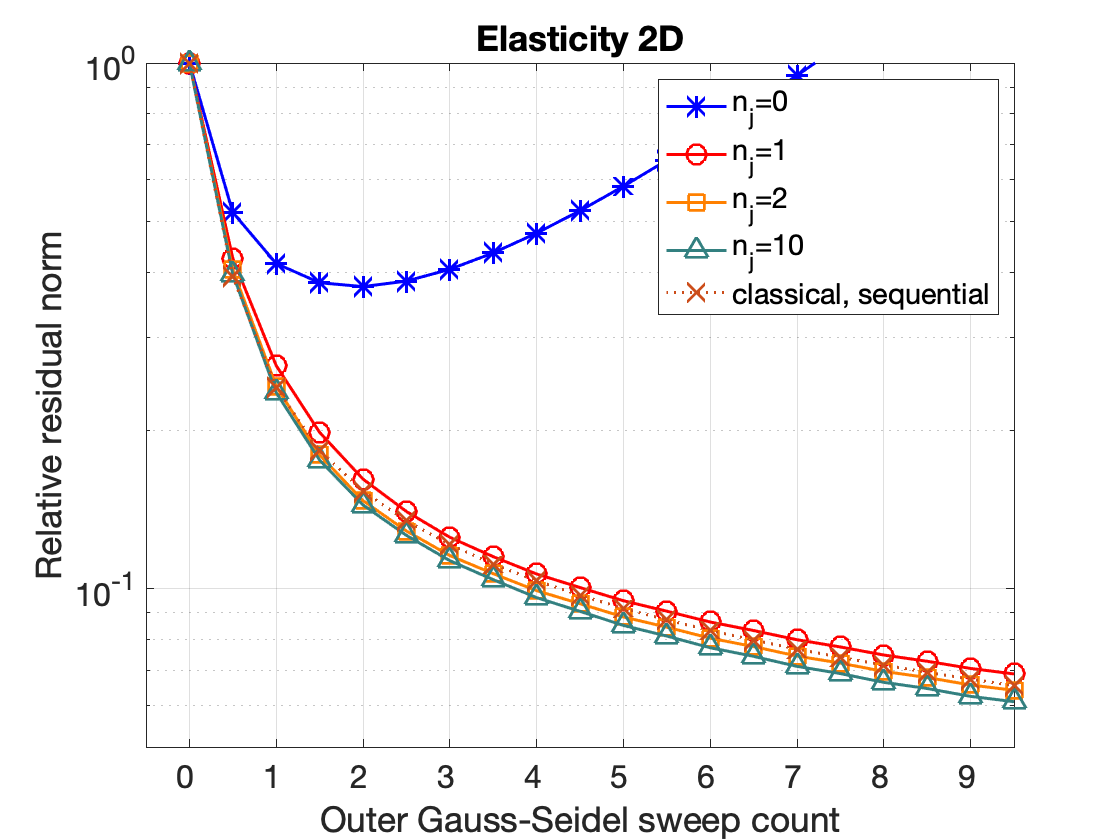}
}
\centerline{
  \includegraphics[width=.5\textwidth]{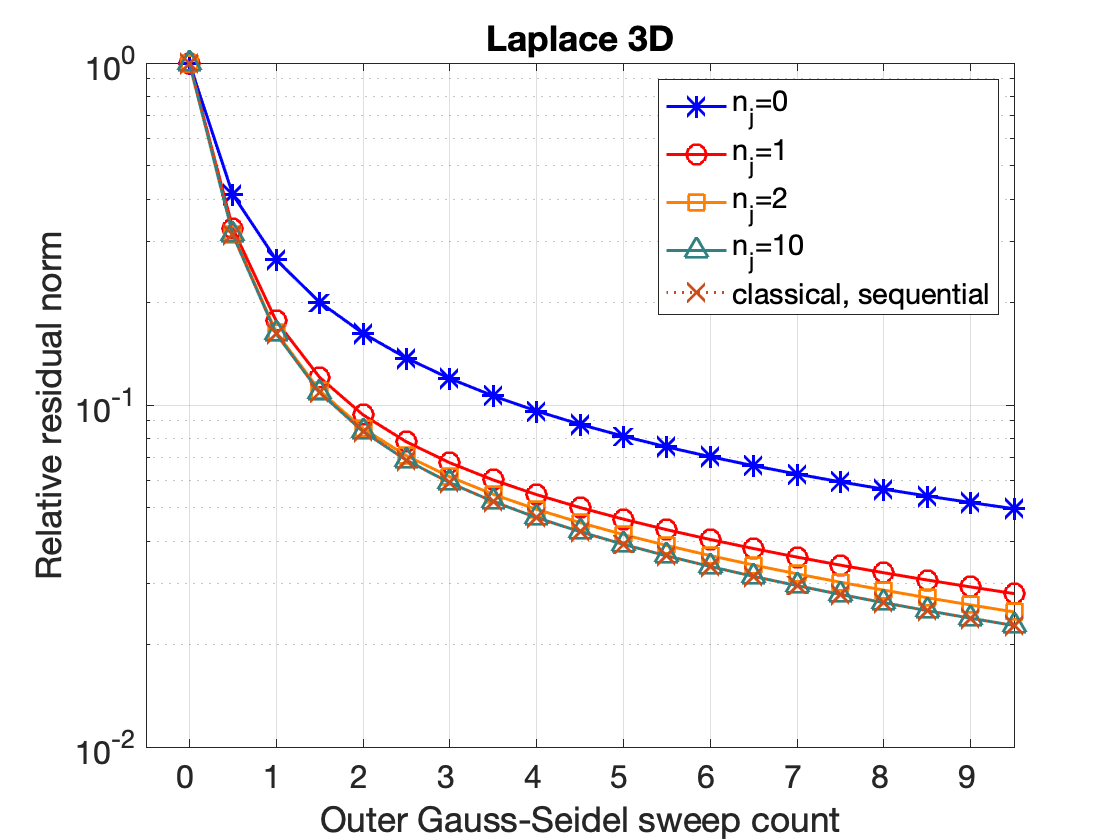}
  \includegraphics[width=.5\textwidth]{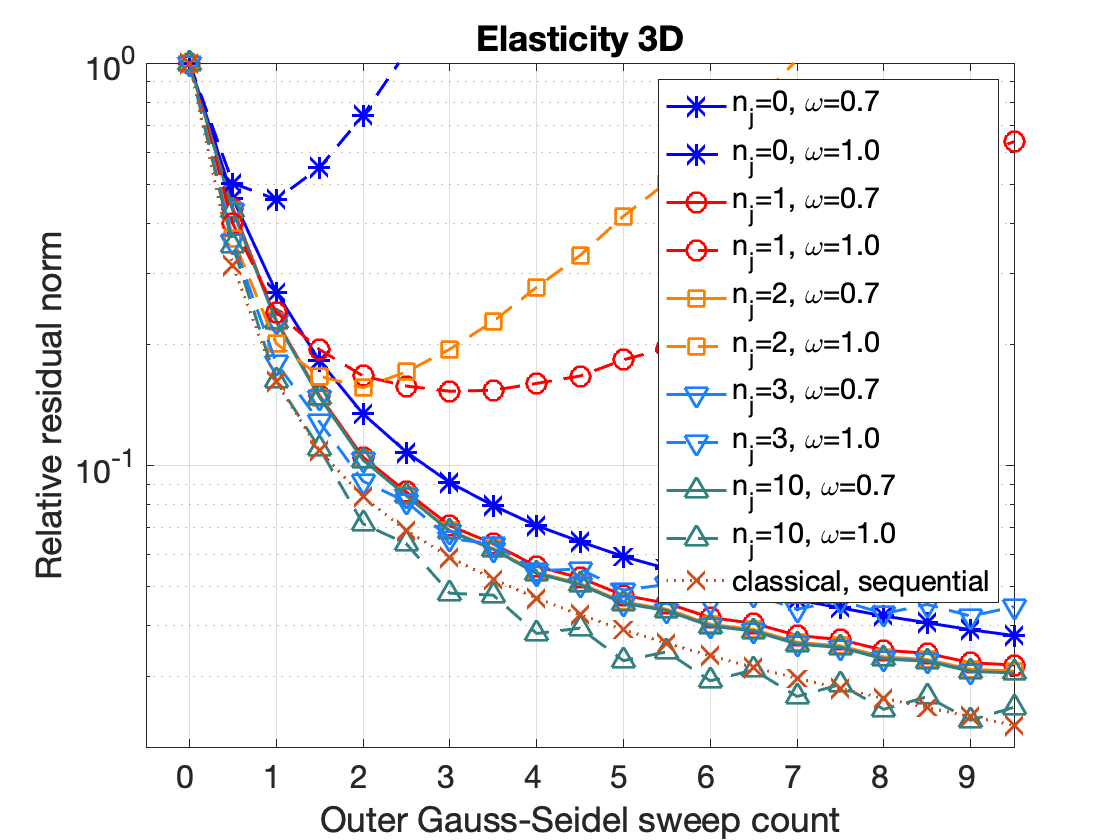}
}
\caption{\label{fig:outer-convergence-model}
Convergence history of $\|\vect{b} - A\vect{x}_k\|_2/\|\vect{b}\|_2$  for SGS2. Elasticity problems.}
\end{figure}
\begin{figure}
\centerline{
  \includegraphics[width=.5\textwidth]{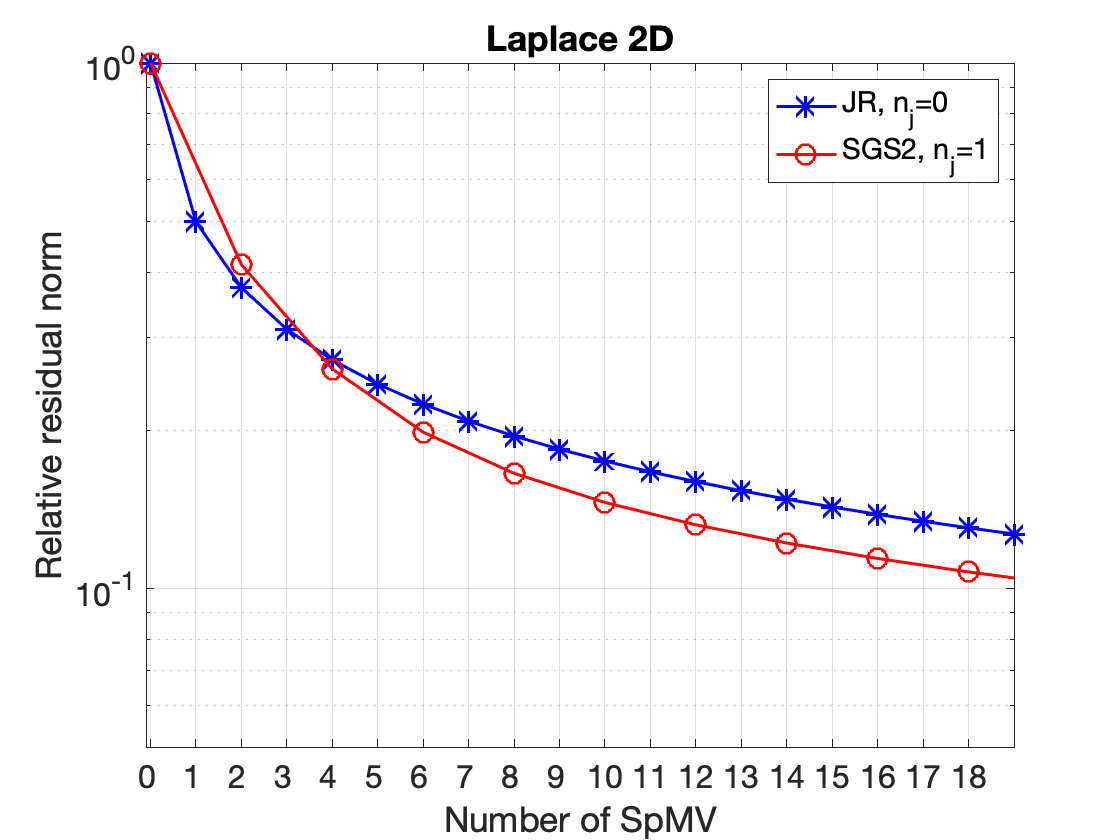}
  \includegraphics[width=.5\textwidth]{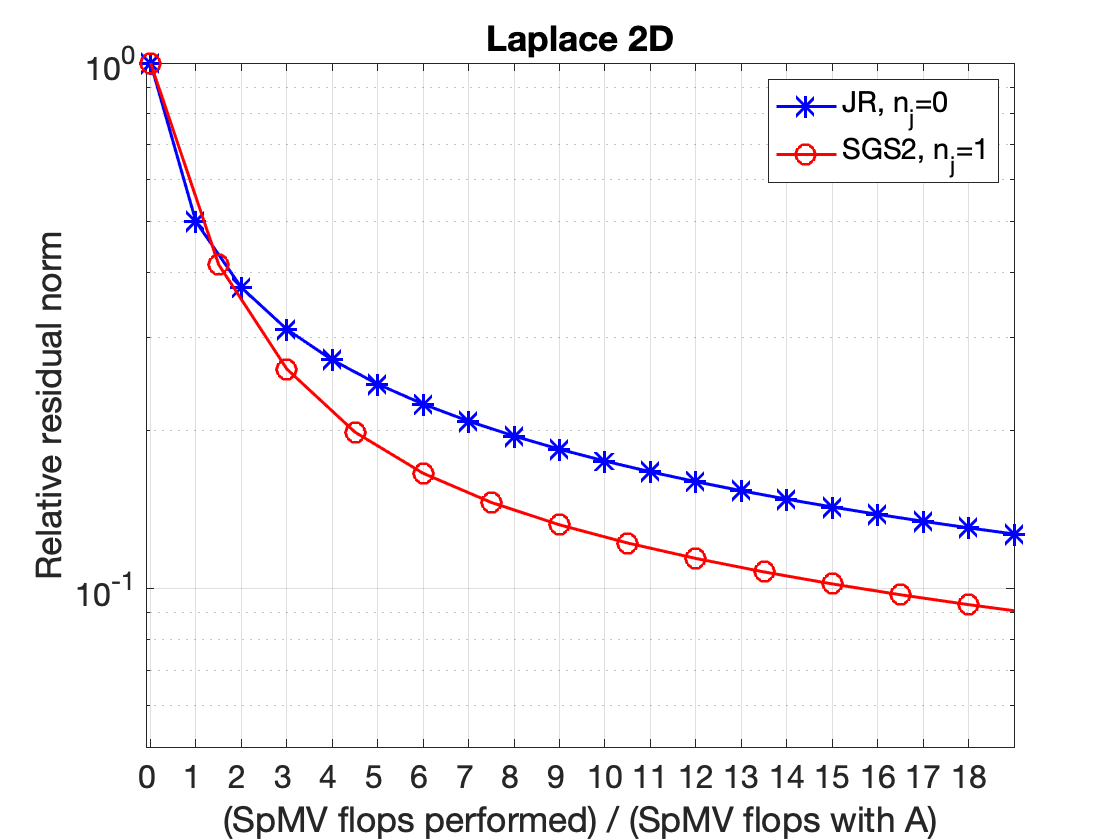}
}
\caption{\label{fig:outer-convergence-model-v2}
Convergence history of $\|\vect{b} - A\vect{x}_k\|_2/\|\vect{b}\|_2$ for GS2
and JR. Laplace 2D problem. Number of SpMV. Ratio of SpMV flops
versus SpMV flops using $A$.}
\end{figure}

First in Figure~\ref{fig:outer-convergence-model}, the convergence of
symmetric SGS2 is displayed as the stand-alone fixed-point iterative
solver, where $n_j = 0$ indicates no inner sweeps, and thus becomes the JR
iterative solver. Our observations are summarized below:
\begin{itemize}
\item
SGS2(1, 1) performs a total of four SpMVs.  
According to the discussion in Section~\ref{sec:others}, 
SGS2(1,1) is likely to converge faster than JR(3), but
it may converge slower than JR(4).  
The difference in the convergence rates was small, but
JR(4) converged faster than SGS2(1, 1) in the first few iterations, but when
more sweeps are performed, SGS2(1, 1) often converged faster than either
JR(3) or JR(4) (as plotted in Figure~\ref{fig:outer-convergence-model-v2}).
\item
JR with the default damping parameter
$\omega=1.0$ diverged, while SGS2 improved the convergence with 
additional inner sweeps.  
With a smaller damping factor $\omega=0.7$, both JR(4) and JR(3) converged faster. 
\item
The figure displays the convergence rate of SGS2 in compact form.  
This requires fewer flops than the non-compact form, 
but may require extra inner sweeps to achieve the same
rate of convergence.  The plot confirms that the compact form
achieves a smaller relative residual norm,  but then stagnates after 
additional outer sweeps.
\item
A small number of inner sweeps ($n_j = 1 \sim 3$) were sufficient for SGS2 to
match the convergence rate of the sequential SGS.
\end{itemize}

\begin{figure}
\centerline{
  \includegraphics[width=.5\textwidth]{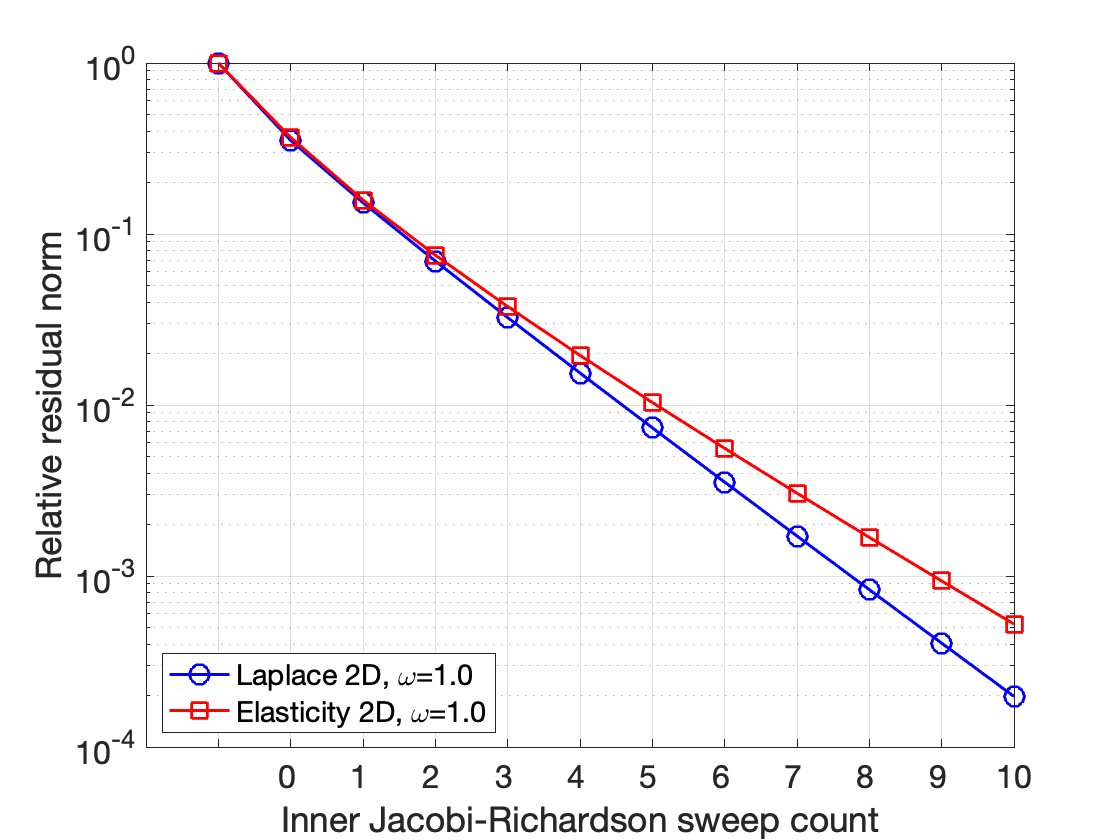}
  \includegraphics[width=.5\textwidth]{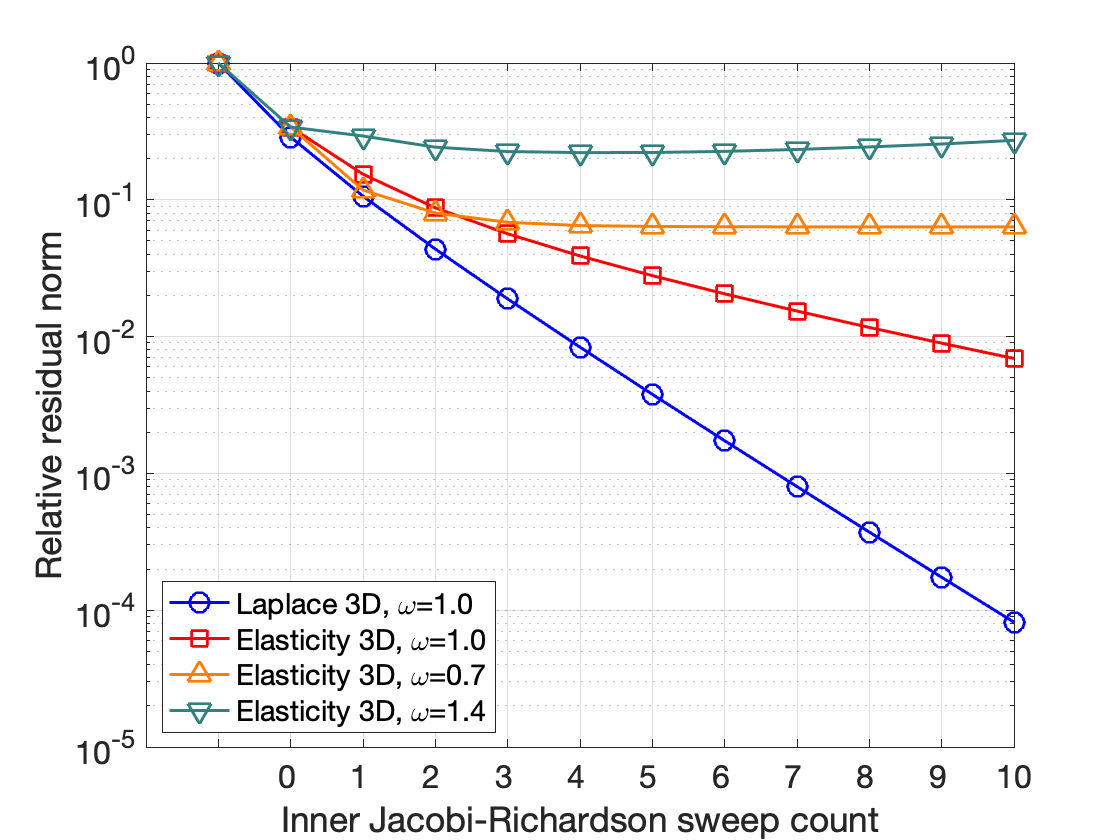}}
\caption{\label{fig:inner-convergence-model}
Convergence history of $\|\vect{r}_l - (L+D)\vect{g}_j\|_2/\|\vect{r}_k\|_2$, 
with inner JR sweeps $(n_t=5, n_j=2)$ for model problems. $n_x=20$.}
\end{figure}

\begin{figure}[htb]
\centerline{
\subfloat[Laplace 2D ($n_x = 1000$)]{
  \includegraphics[width=.5\textwidth]{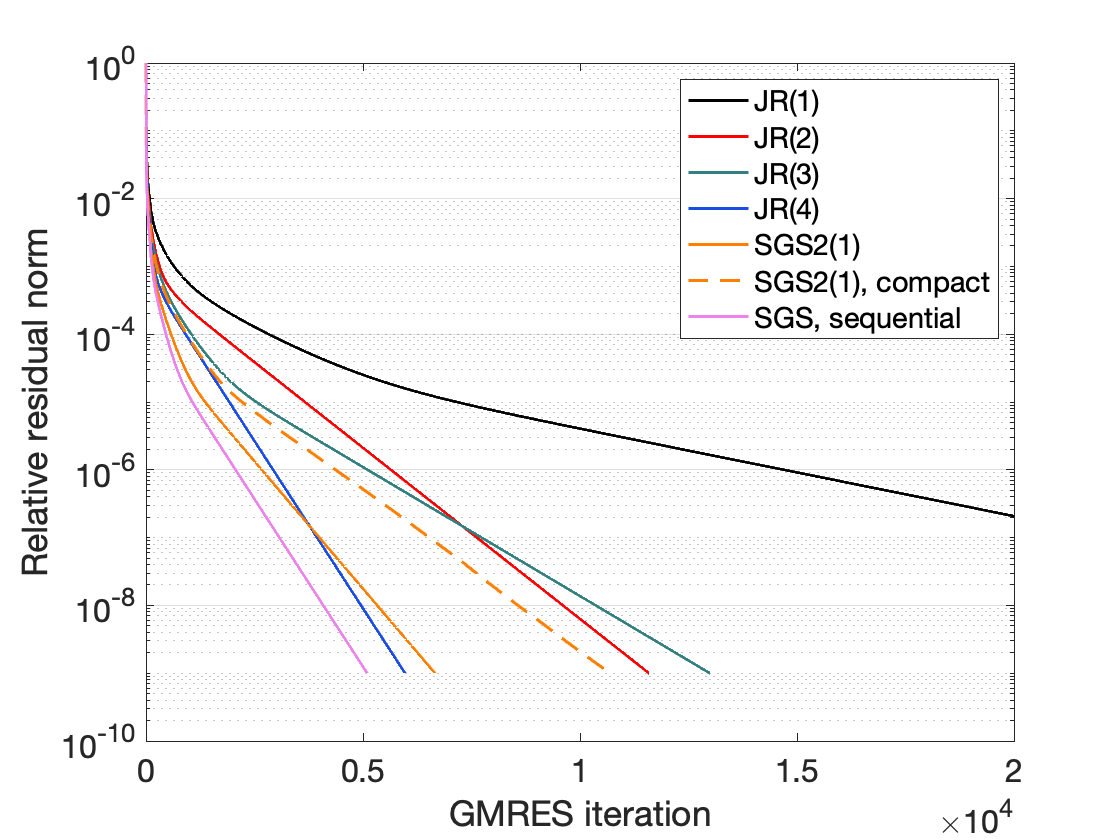}}
\subfloat[Laplace 3D ($n_x = 100$)]{
  \includegraphics[width=.5\textwidth]{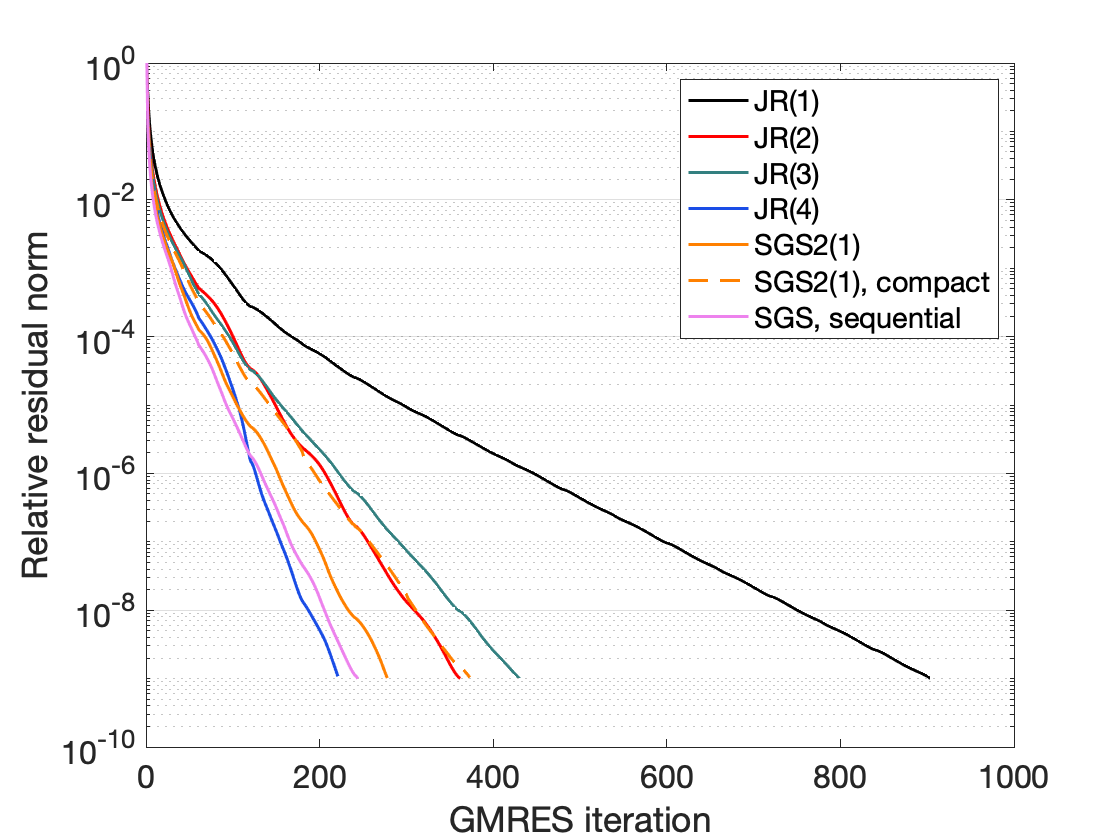}}
}
\caption{\label{fig:gmres-convergence-model}
GMRES(60) convergence with Jacobi Richardson preconditioner.}
\end{figure}

Figure~\ref{fig:inner-convergence-model}  displays the inner JR 
convergence rate for solving the triangular system, while
Figure~\ref{fig:gmres-convergence-model} displays  the
convergence history with the SGS2 preconditioner.
One to three JR sweeps were enough to 
reduce the relative residual norm by an order of magnitude, which is often 
enough to obtain the desired convergence rate.
The iteration count is reduced by increasing the number of
JR sweeps. The SGS2(1,1) iteration results in a faster 
convergence rate than JR(3) while it leads to
a similar  rate as JR(4), or sequential SGS(1).
The latter degraded using the SGS2 preconditioner in compact form.

\begin{figure}[htb]
\centerline{
\subfloat[Laplace 2D]{
  \includegraphics[width=.5\textwidth]{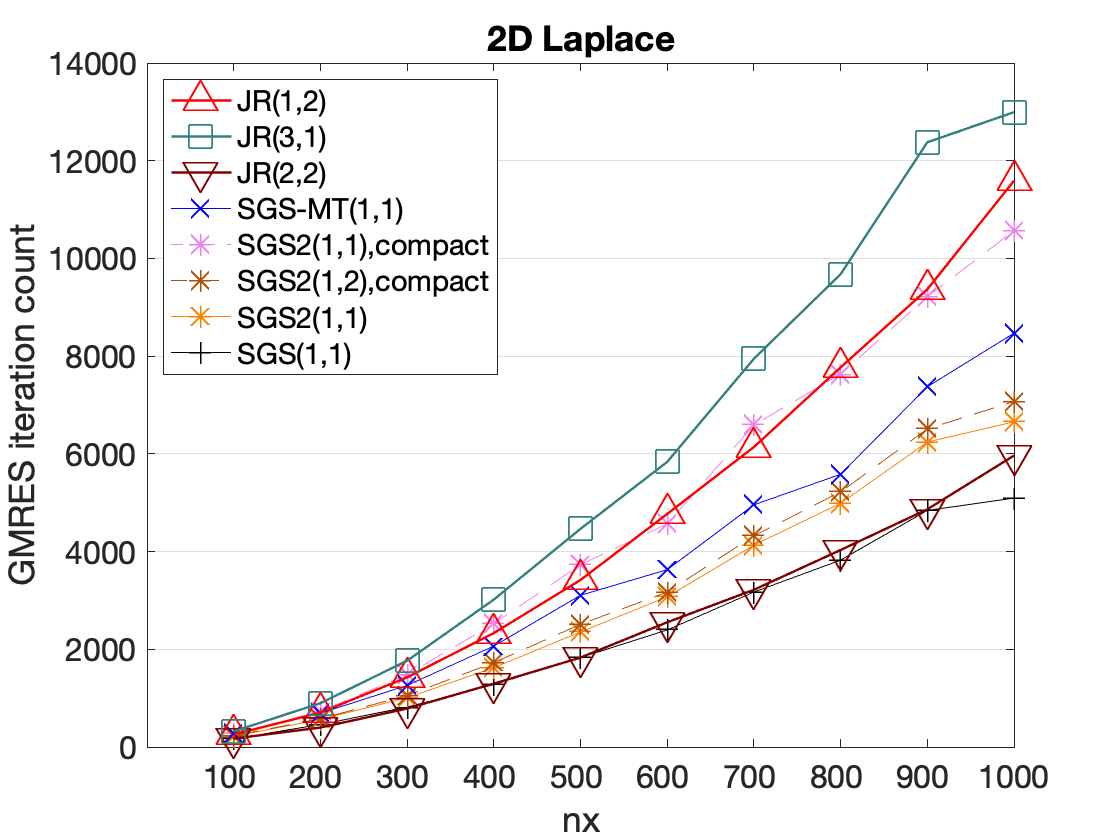}}
\subfloat[Laplace 3D]{
  \includegraphics[width=.5\textwidth]{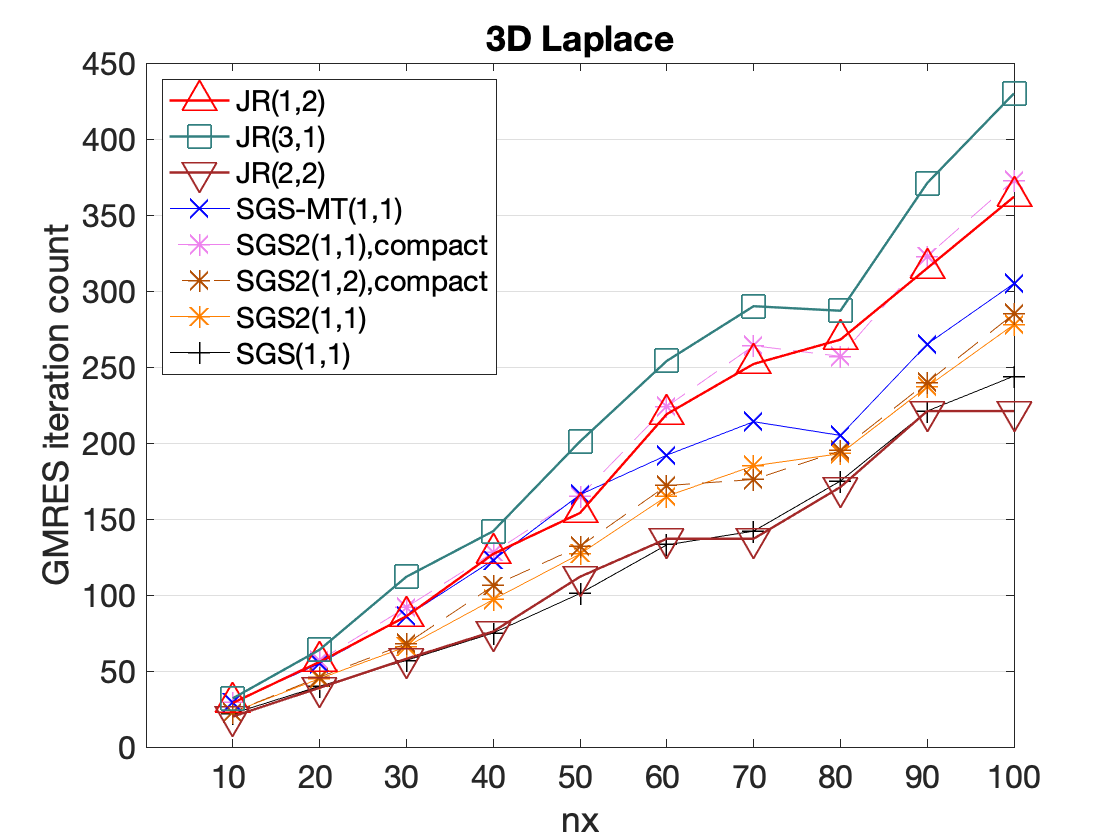}}
}
\caption{\label{fig:laplace-iter}
GMRES(60) iteration counts for Laplace problems.}
\end{figure}

\begin{figure}[htb]
\centerline{
\subfloat[Laplace 2D]{
  \includegraphics[width=.5\textwidth]{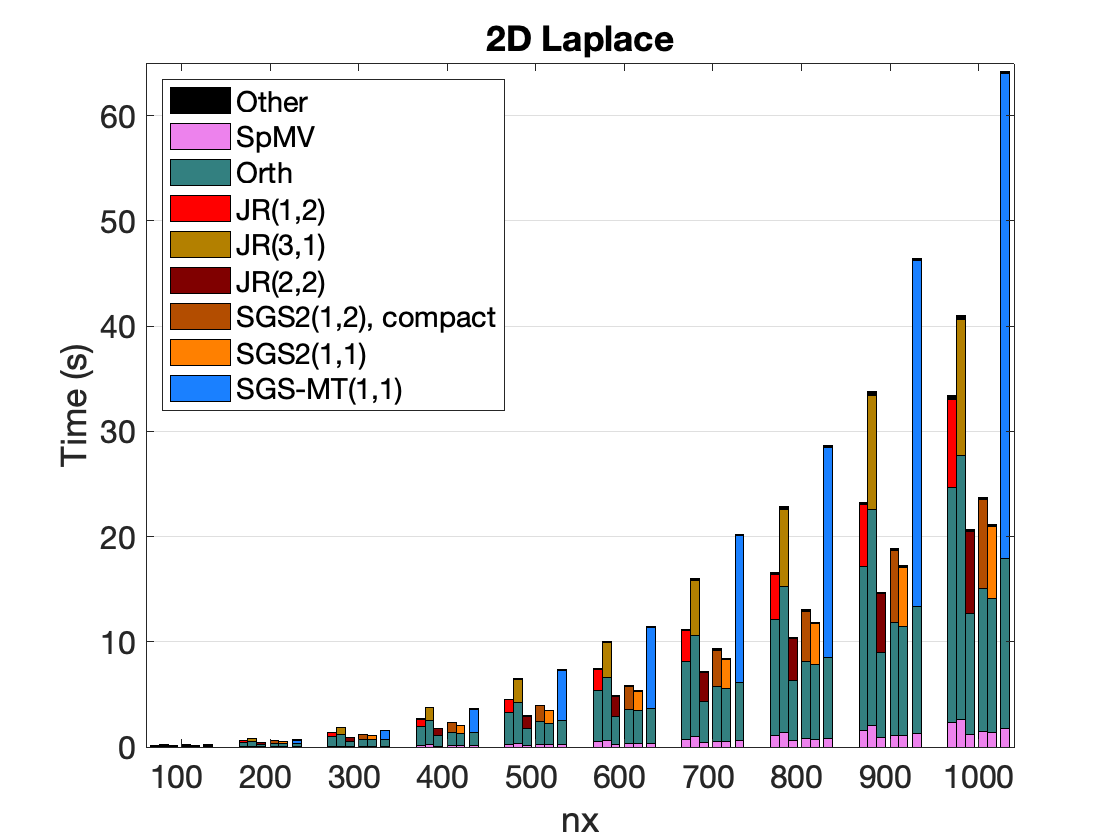}}
\subfloat[Laplace 3D]{
  \includegraphics[width=.5\textwidth]{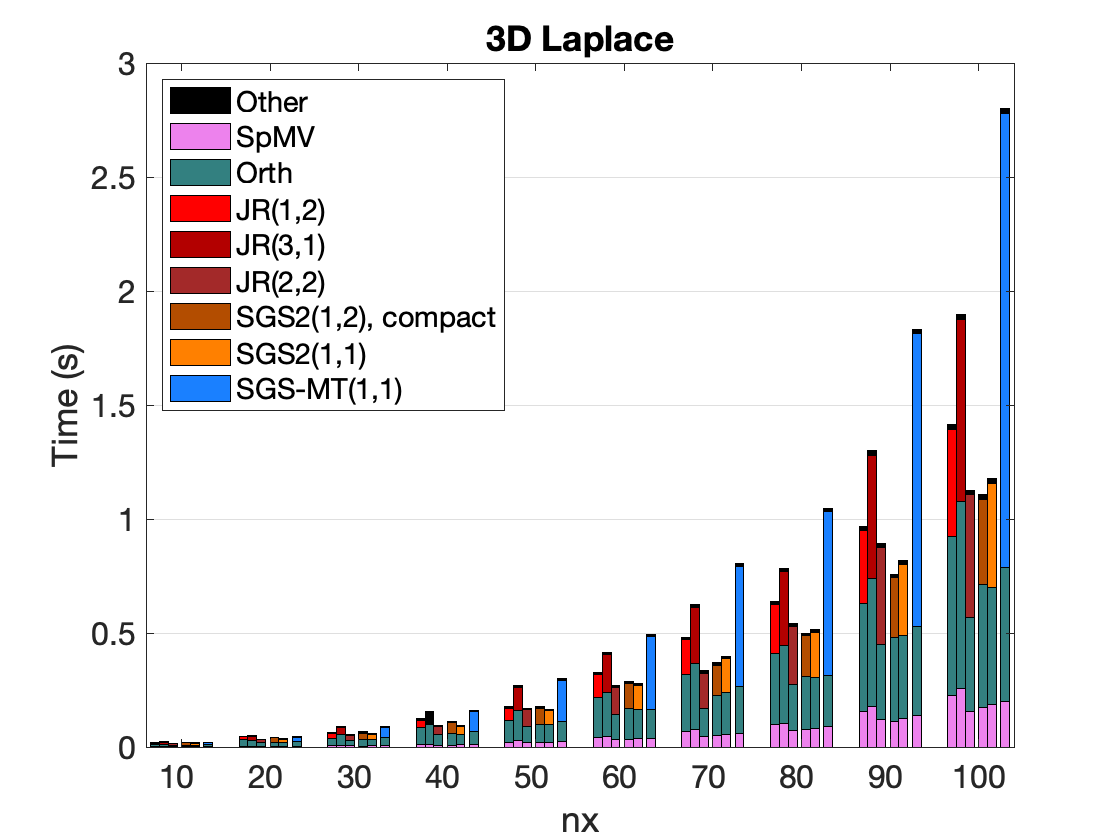}}
}
\caption{\label{fig:laplace-time}
GMRES(60) solution time on one GPU for Laplace problems.}
\end{figure}

\begin{figure}[htb]
\centerline{
\subfloat[Laplace 2D]{
  \includegraphics[width=.5\textwidth]{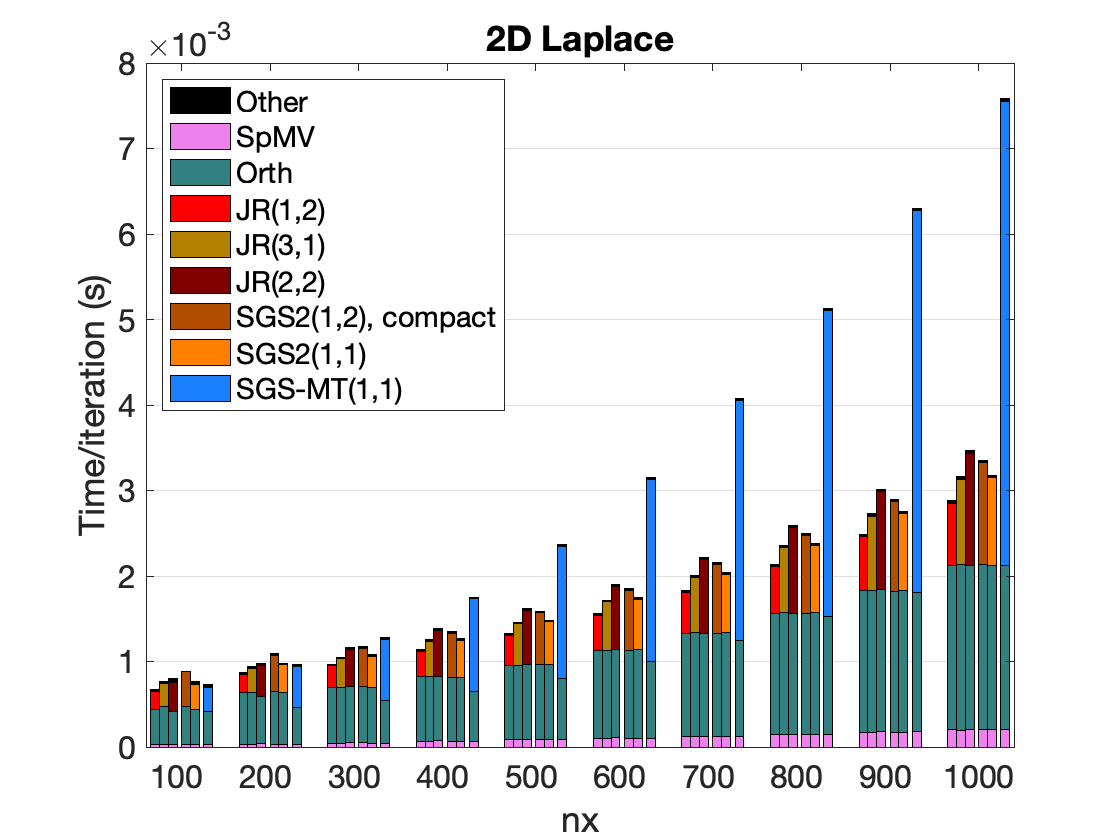}}
\subfloat[Laplace 3D]{
  \includegraphics[width=.5\textwidth]{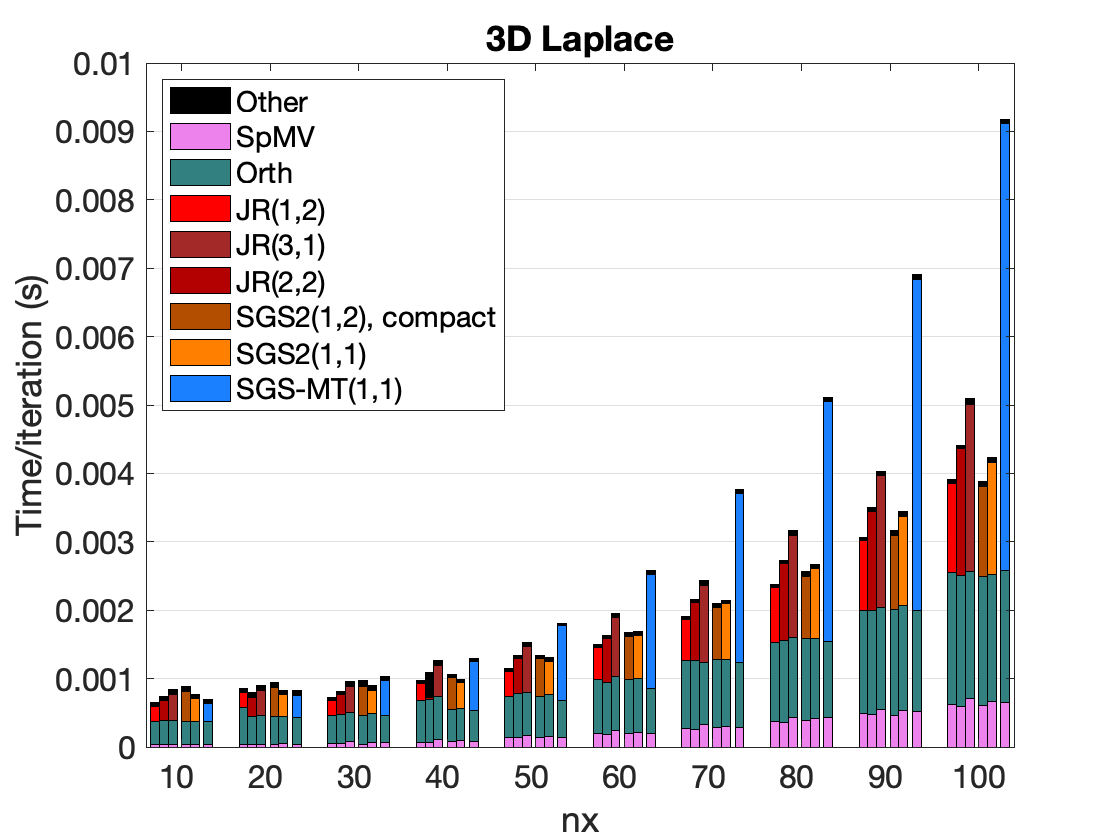}}
}
\caption{\label{fig:laplace-iter-time}
GMRES(60) time per iteration on one GPU for Laplace problems.}
\end{figure}

The time-to-solution with the SGS2 preconditioner was evaluated 
for model problems. Figures~\ref{fig:laplace-iter} through
\ref{fig:laplace-iter-time} compare the iteration count, time to
solution, and time per iteration for the Laplace problems.  
Although a variable number of iterations is required to converge,
the computed times of the different solvers were similar. 
To summarize,
\begin{itemize} 
\item
The solver converged with fewer iterations using SGS2(1,1)
versus JR(3), but required more iterations than with JR(4),
which achieved a convergence rate similar to the sequential SGS(1).
\item
In many cases, SGS2 requires more iterations than the
sequential SGS, but fewer iterations than the MT-SGS, which
could increase the iteration count.
\item
Compared with SGS2(1, 1), JR(4) exhibited a larger compute time
per iteration but needed fewer iterations to converge.  Overall,
SGS2(1,1) and JR(4) obtained similar time to solution.  
\item
Because SGS2 only needed a small number of SpMVs to obtain convergence
rates similar to the sequential SGS, SGS2 was faster than sequential SGS
or multi-threaded MT-SGS.  Furthermore, an SpMV is inherently parallel,
and thus SGS2 may also exhibit a slower increase in compute time with
the problem size..
\end{itemize}


\begin{figure}[htb]
\centerline{
\subfloat[2D problems.\label{fig:2d-elasticity-iters}]{
  \includegraphics[width=.5\textwidth]{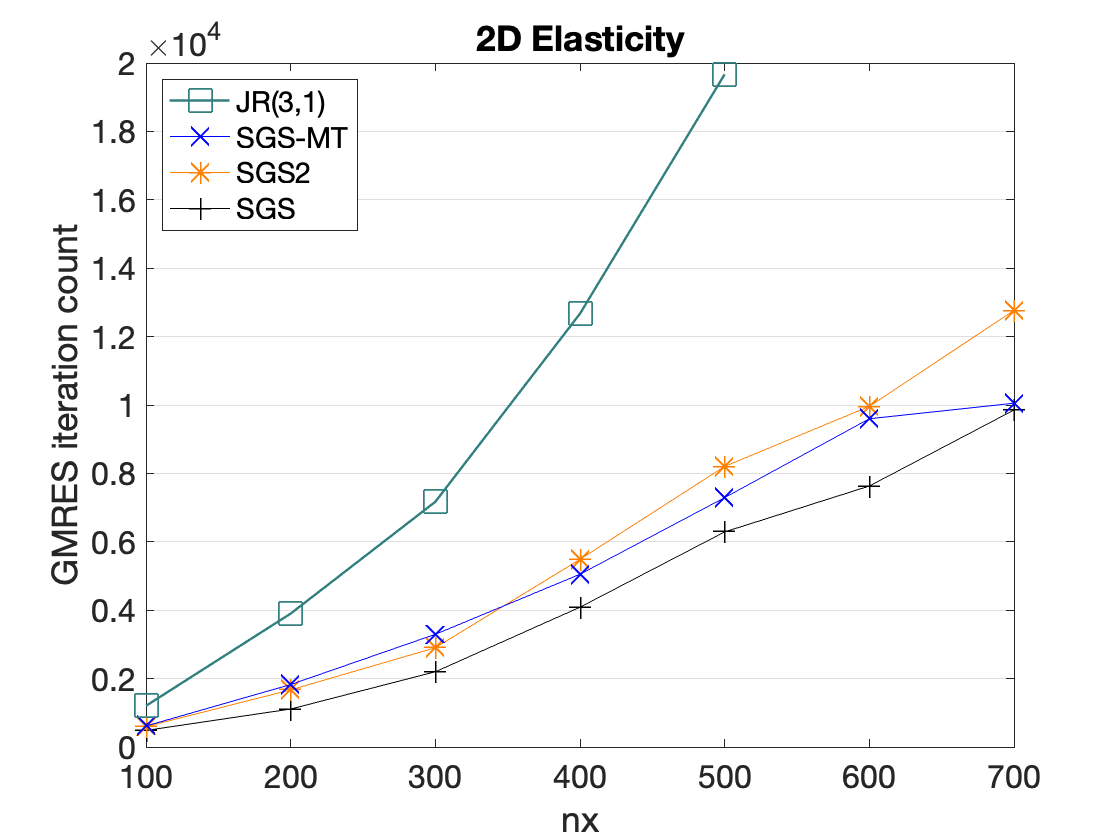}}
\subfloat[3D problems.\label{fig:3d-elasticity-iters}]{
  \includegraphics[width=.5\textwidth]{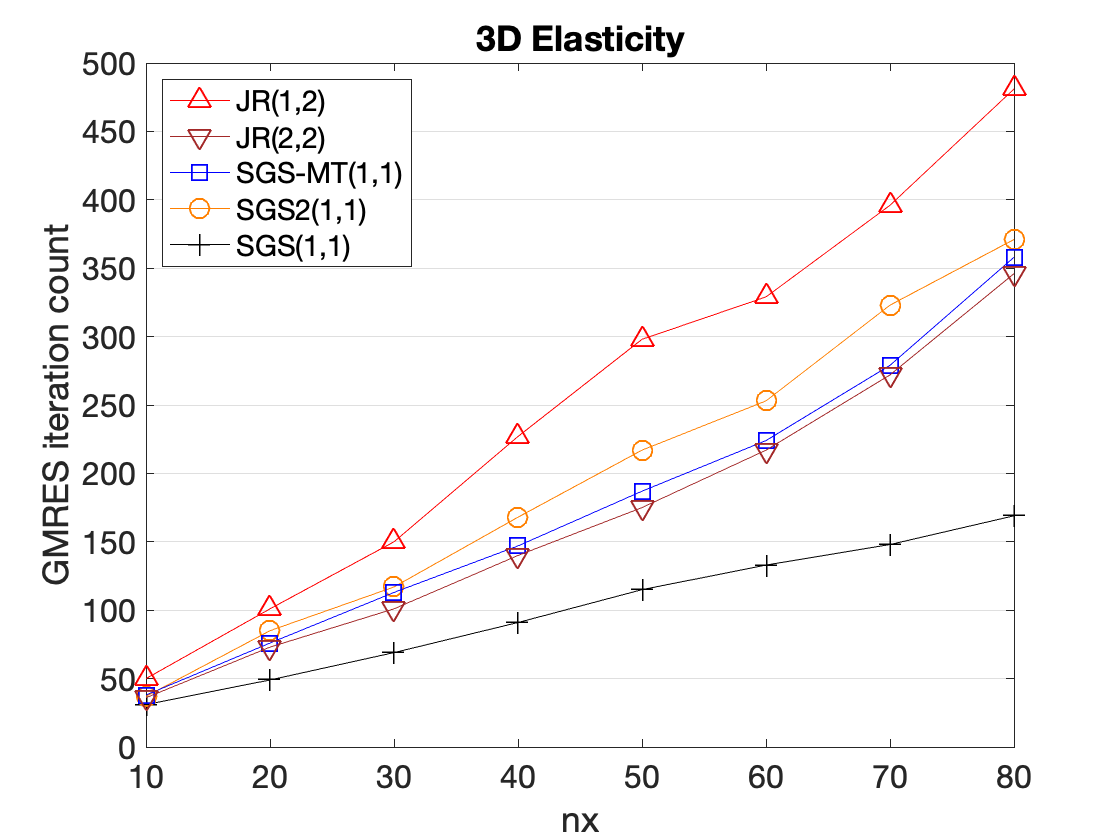}}
}
\caption{GMRES iterations for elasticity problems.}
\end{figure}

\begin{figure}[htb]
\centerline{
\subfloat[GMRES time to solution]{
  \includegraphics[width=.5\textwidth]{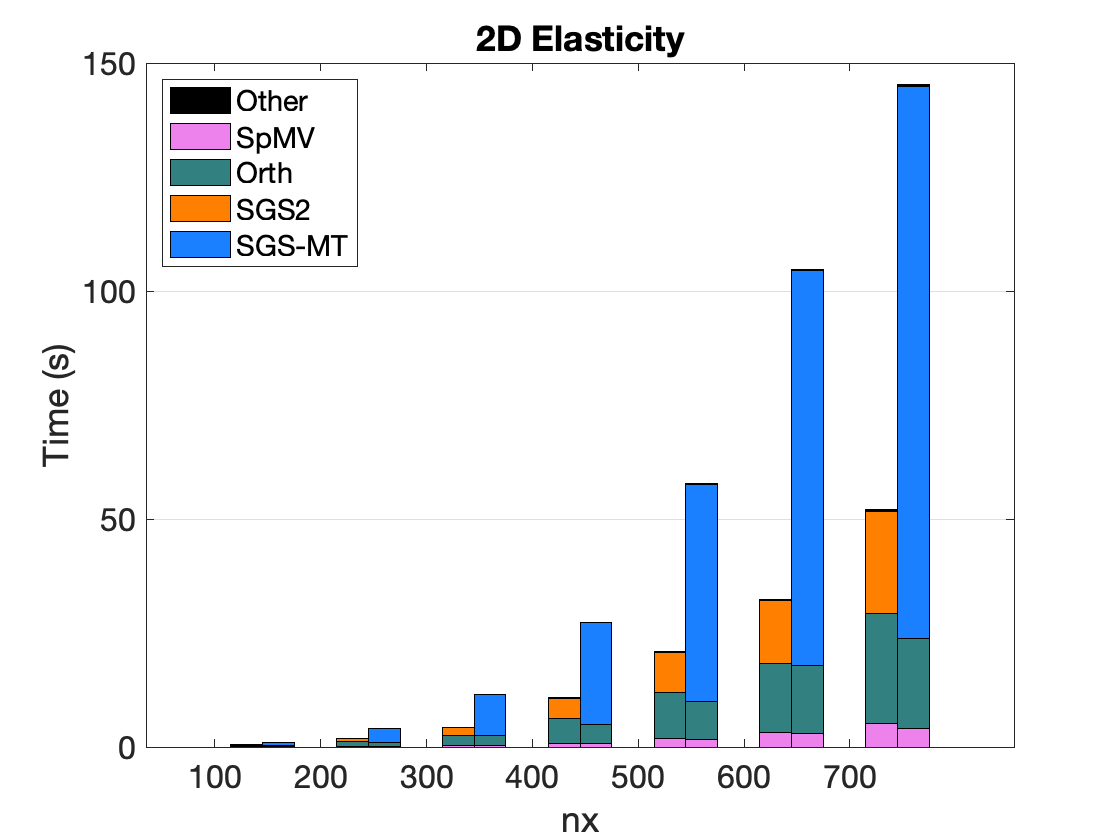}}
%
\subfloat[GMRES time per iteration]{
  \includegraphics[width=.5\textwidth]{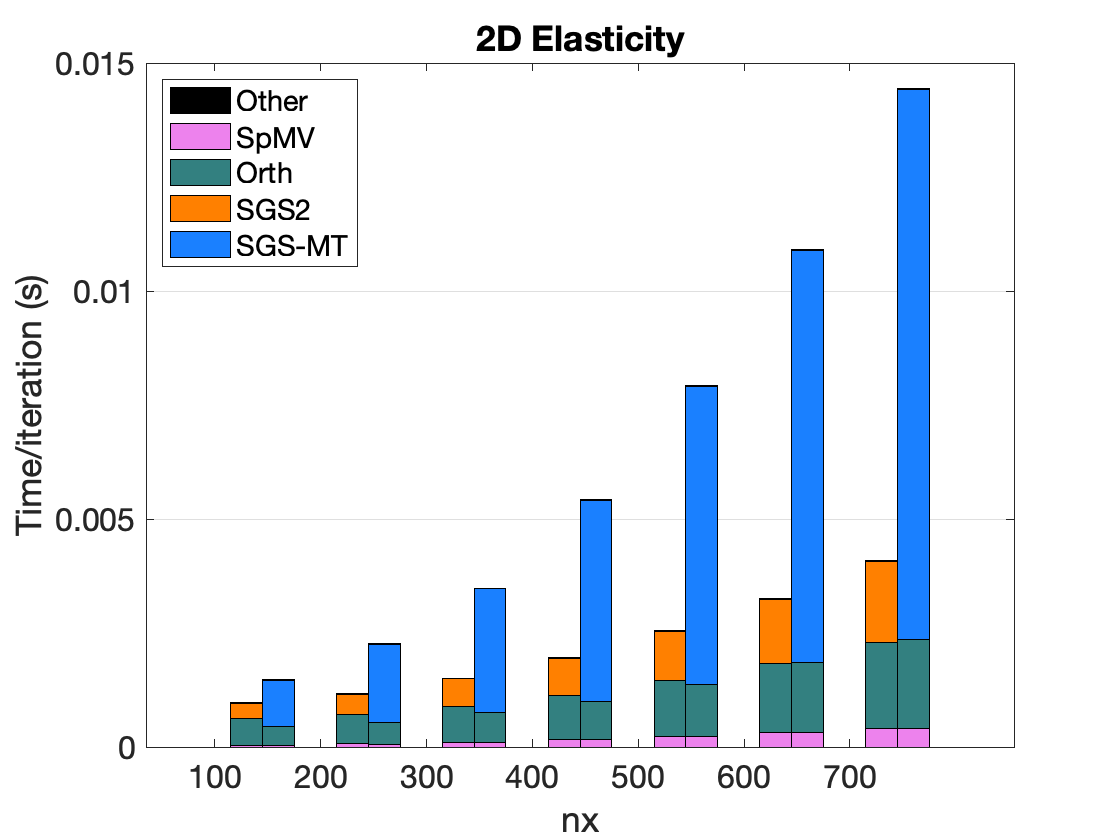}}
}
\caption{Compute times with one GPU for 2D elasticity problems.}
\label{fig:2d-elasticity}
\end{figure}

Figures~\ref{fig:2d-elasticity-iters} and~\ref{fig:2d-elasticity} display
iteration counts and compute times for the SGS2 preconditioner when applied to
elasticity problems. The solver failed to converge using either JR(1,2) or JR(2,2),
with the default damping factor of $\omega=1.0$.
For the 3D problems, the damping factors were carefully selected to be
(1.4, 1.2, 0.8, 0.7). To compare, the sequential SGS, MT-SGS, and
SGS2 required about the same number of iterations, but SGS2 results 
in a much shorter time to solution.


\begin{figure}[htb]
\centerline{
\subfloat[GMRES solution time and time per iteration]{
  \includegraphics[width=.5\textwidth]{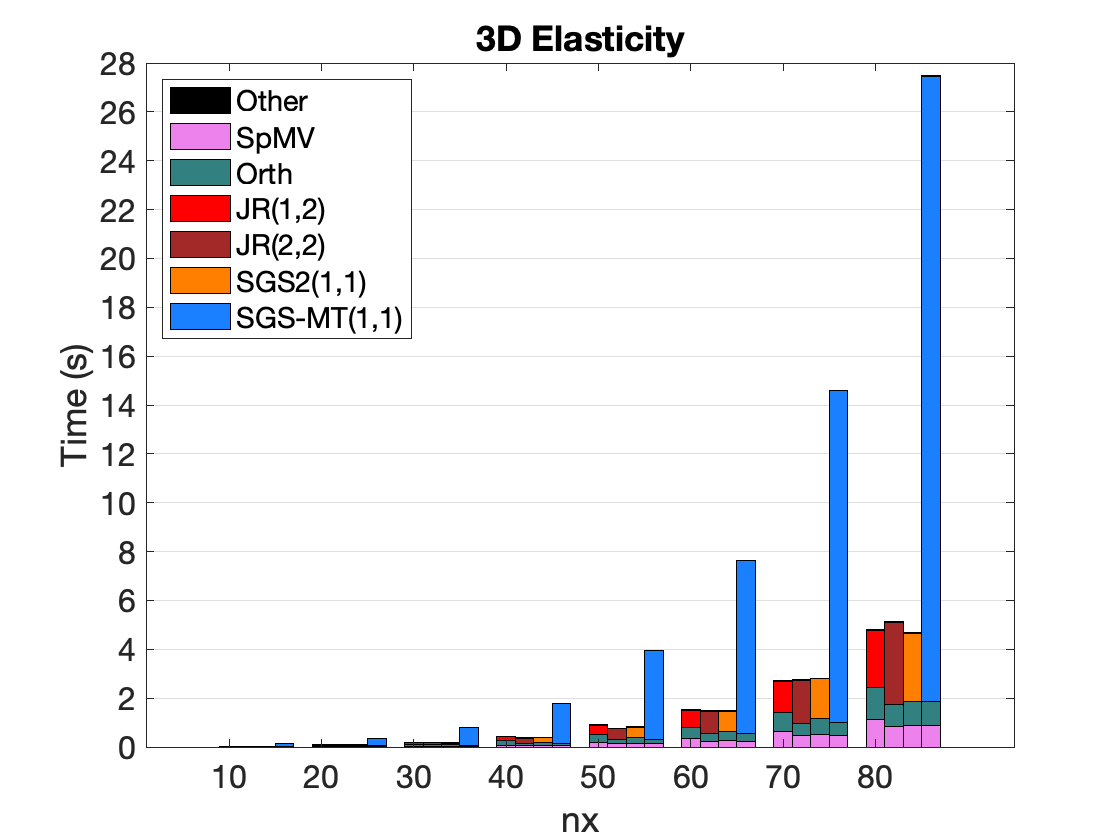}
  \includegraphics[width=.5\textwidth]{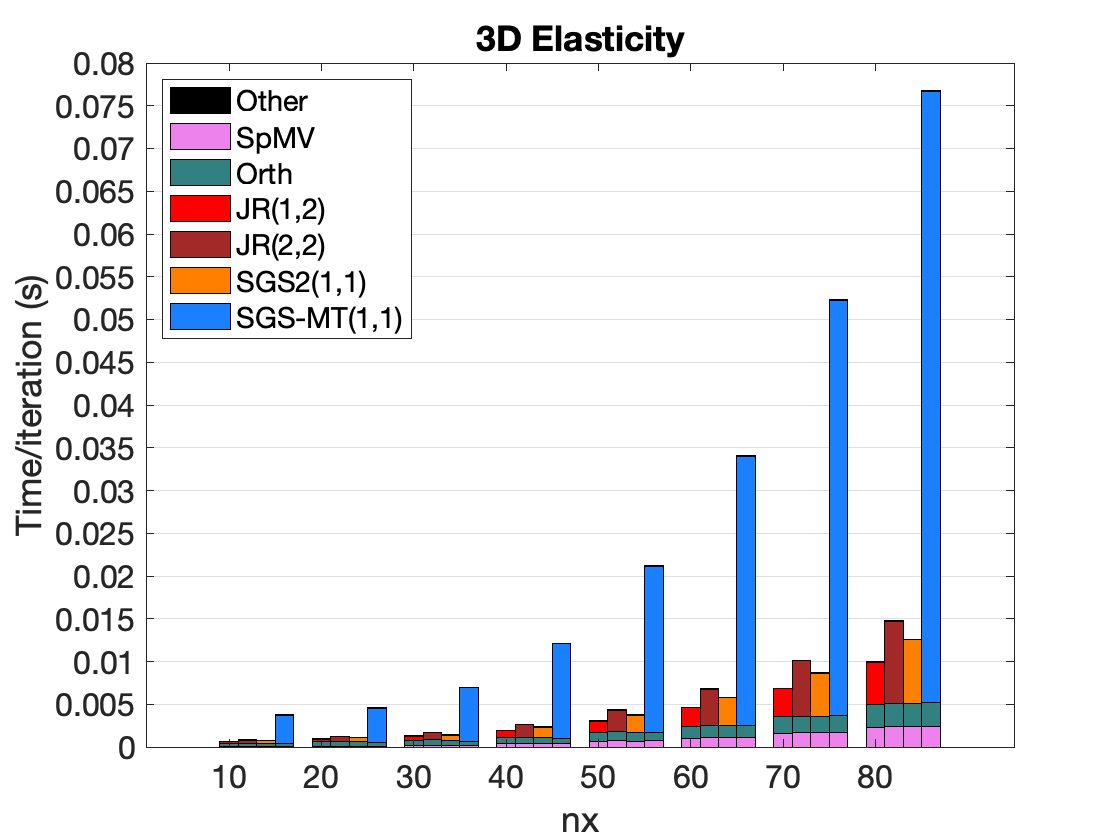}}
}
\caption{Solution time with one GPU for 3D elasticity problems.}
\label{fig:3d-elasticity}
\end{figure}

Next, the convergence rate with the SGS2 preconditioner applied to
elasticity problems was studied, with $n_x=500$ and $n_x=20$.
These were evaluated by varying the number of inner sweeps (i.e. $n_j$).
SGS2 required a few inner sweeps to converge when a smaller damping factor 
is employed (see Figure~\ref{fig:outer-convergence-model}). 
Figure~\ref{fig:3d-elasticity-iters} displays the iteration counts and
Figure~\ref{fig:3d-elasticity} provides the compute times for 3D 
elasticity problems.

Figure~\ref{tab:ela3d-omega}
displays the iteration counts with different values of the damping factor
and a varying number of inner or outer sweeps for the SGS2 and JR,
respectively.  Only one outer sweep was performed for
SGS, MT-SGS, and SGS2. In summary,

\begin{itemize}
\item
With the default damping factor $\omega=1$, when SGS2
performs enough inner iterations, the same convergence rate is achieved as
the sequential SGS.  
\item
When $\omega=1$, JR diverges for this problem. As a result, even with multiple JR
sweeps, the solver convergence was erratic.
\item
Conversely, JR converges well with a smaller damping factor $\omega=0.7$. 
Using a carefully-selected damping factor, the 
SGS2(1,1) iteration count was less than JR(3) but larger than JR(4).
\item 
With increasing inner sweeps for SGS2, the iteration count is the same 
as sequential SGS.
\end{itemize}
%

\begin{figure}
\begin{center}\footnotesize
\begin{tabular}{l|rrrrrrrrrrrr}
$\omega$   & 1.6 & 1.5 & 1.4      & 1.3 & 1.2       & 1.1 & 1.0      & 0.9      & 0.8      & 0.7       & 0.6       & 0.5\\ 
\hline\hline
SGS        & 52  & 49  & {\bf 49} & 51  & 53        & 56  & 61       & 68       & 78       & 87        &  93       & 102\\
SGS-MT     & 93  & 90  & 83       & 80  & {\bf 77}  & 78  & 79       & 78       & 80       & 87        &  98       & 104\\
\hline\hline
SGS2(1,1)  & --  & --  & --       & --  & --        & --  & --       & 93       & {\bf 85} & 87        &  96       & 103\\
SGS2(1,2)  & --  & --  & --       & --  & --        & --  & --       & {\bf 77} & 82       & 87        &  94       & 103\\
SGS2(1,3)  & --  & --  & --       & --  & --        & --  & 163      & {\bf 68} & 80       & 87        &  93       & 102\\
SGS2(1,4)  & --  & --  & --       & --  & --        & --  & 159      & {\bf 68} & 79       & 87        &  93       & 102\\
SGS2(1,5)  & --  & --  & --       & --  & --        & 2650& {\bf 62} & 68       & 79       & 87        &  93       & 102\\
SGS2(1,10) & --  & --  & --       & --  & --        & 92  & {\bf 61} & 68       & 78       & 87        &  93       & 102\\
\hline\hline
JR(1)      & 213 & 213 & 213      & 213 & 213       & 213 & 213      & 213      & 213      & 213       & 213       & 213\\
JR(2)      & --  & --  & --       & --  & --        & --  & --       & --       & --       & {\bf 101} & 108       & 119\\
JR(3)      & 457 & 399 & 341      & 265 & 233       & 211 & 194      & 152      & 122      & 107       & {\bf 100} & 101\\
JR(4)      & --  & --  & --       & --  & --        & --  & --       & --       & --       & {\bf 73}  & 76        & 82\\
JR(5)      & 1800& 1332& 957      & 696 & 488       & 347 & 227      & 155      & 104      & 78        & {\bf 73}  & 75\\
JR(6)      & --  & --  & --       & --  & --        & --  & --       & --       & --       & {\bf 56}  & 60        & 67\\
JR(10)     & --  & --  & --       & --  & --        & --  & --       & --       & --       & {\bf 43}  & 46        & 51\\
\end{tabular}
\caption{GMRES iteration counts with one GPU for 3D elasticity
problems, $n_x = 20$.  Outer and inner sweeps ($n_t$ and $n_j$ in Fig.~\ref{code:two-stage}).}\label{tab:ela3d-omega}
\end{center}
\end{figure}

\begin{figure}[htb]
\centerline{
  \includegraphics[width=.5\textwidth]{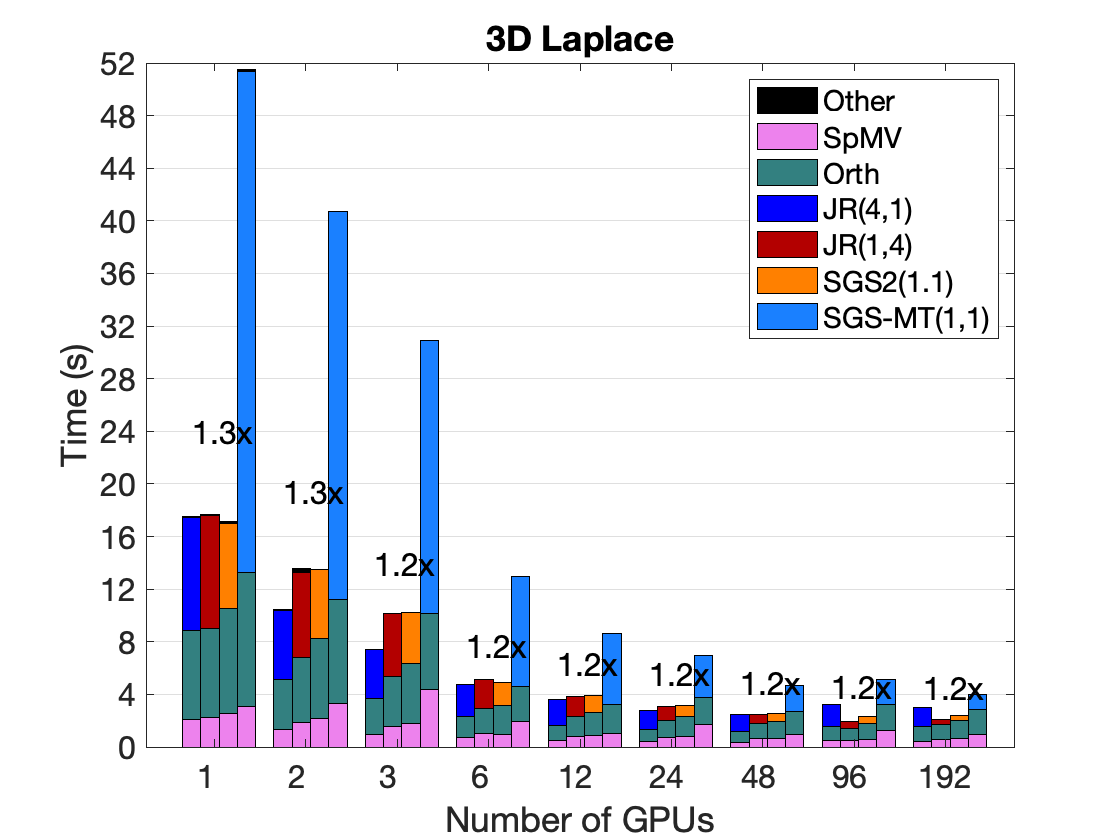}
  \includegraphics[width=.5\textwidth]{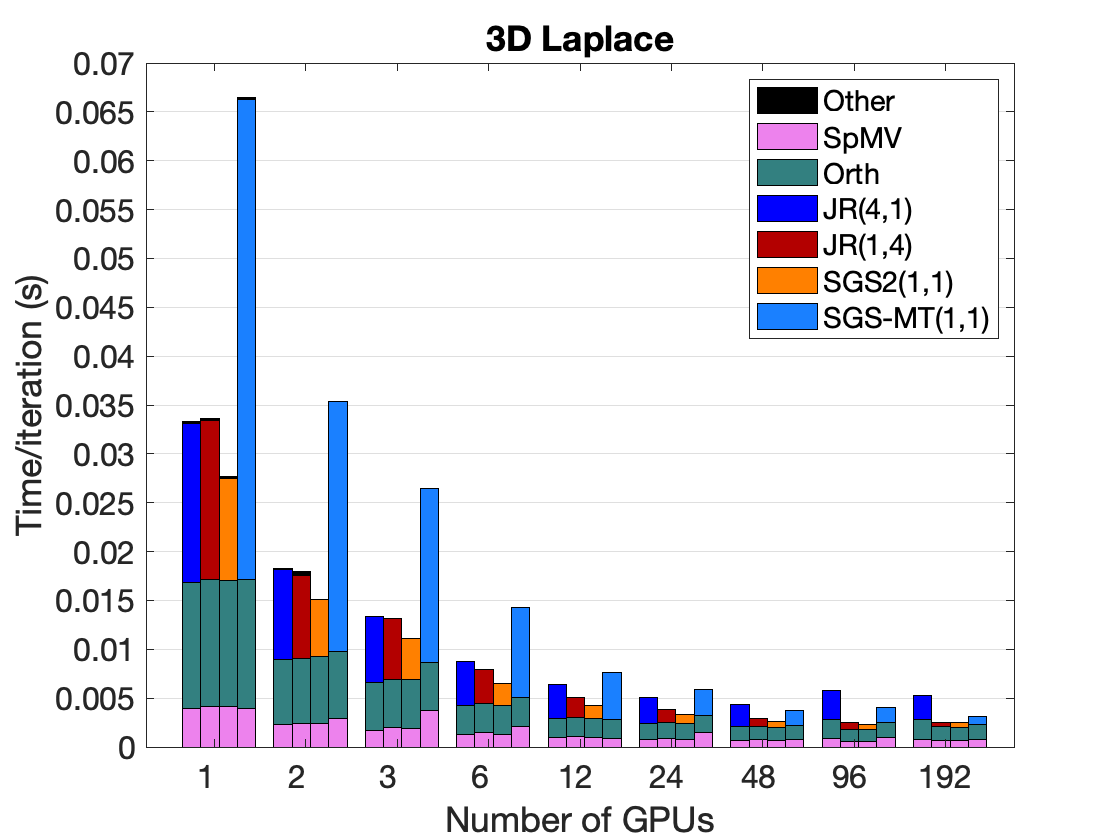}}
\caption{GMRES(60) strong scaling for Laplace 3D problems with $n_x = 200$. 
}
\label{fig:laplace3d-strong-scale}
\end{figure}

\begin{figure}[htb]
\centerline{
  \includegraphics[width=.4\textwidth]{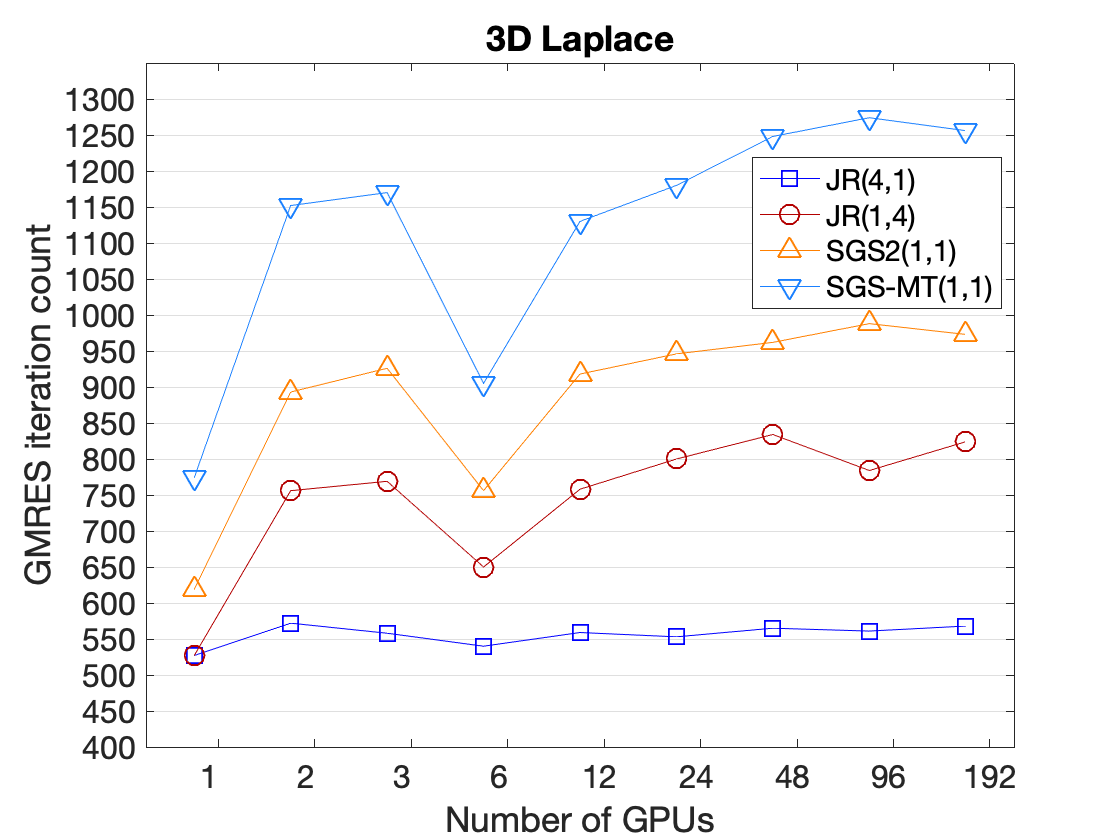}
  \includegraphics[width=.4\textwidth]{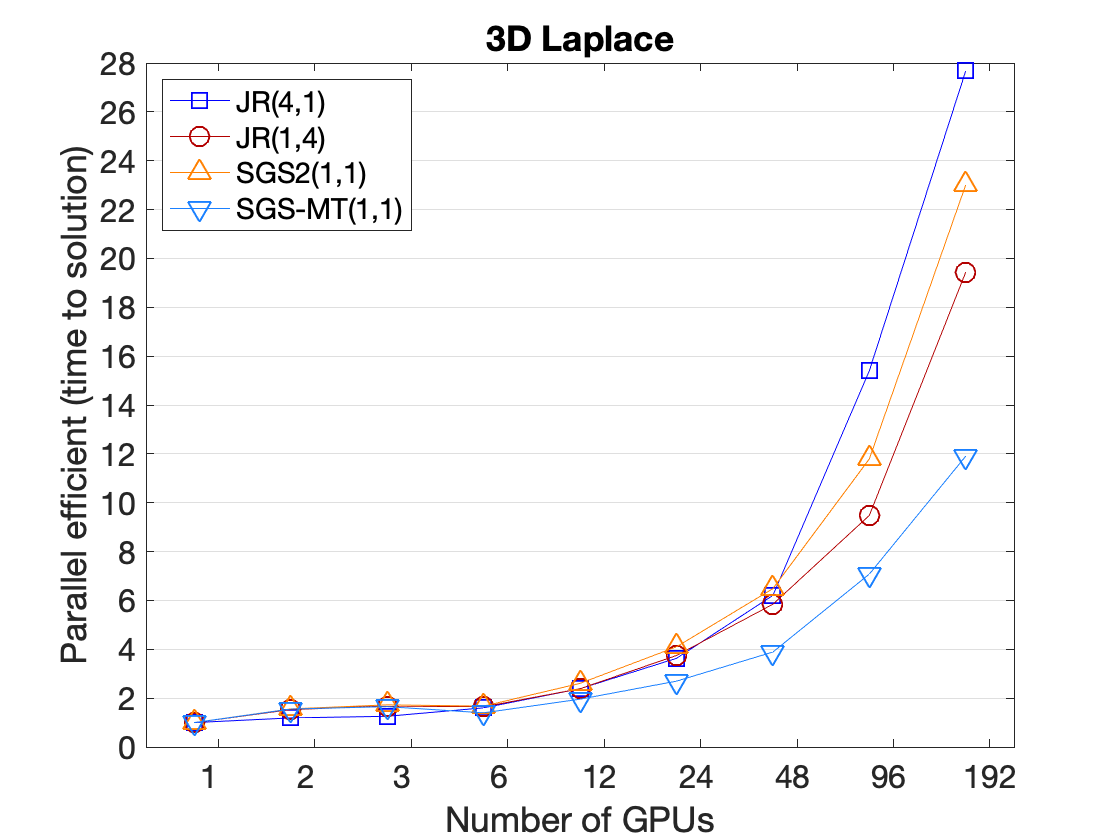}
  \includegraphics[width=.4\textwidth]{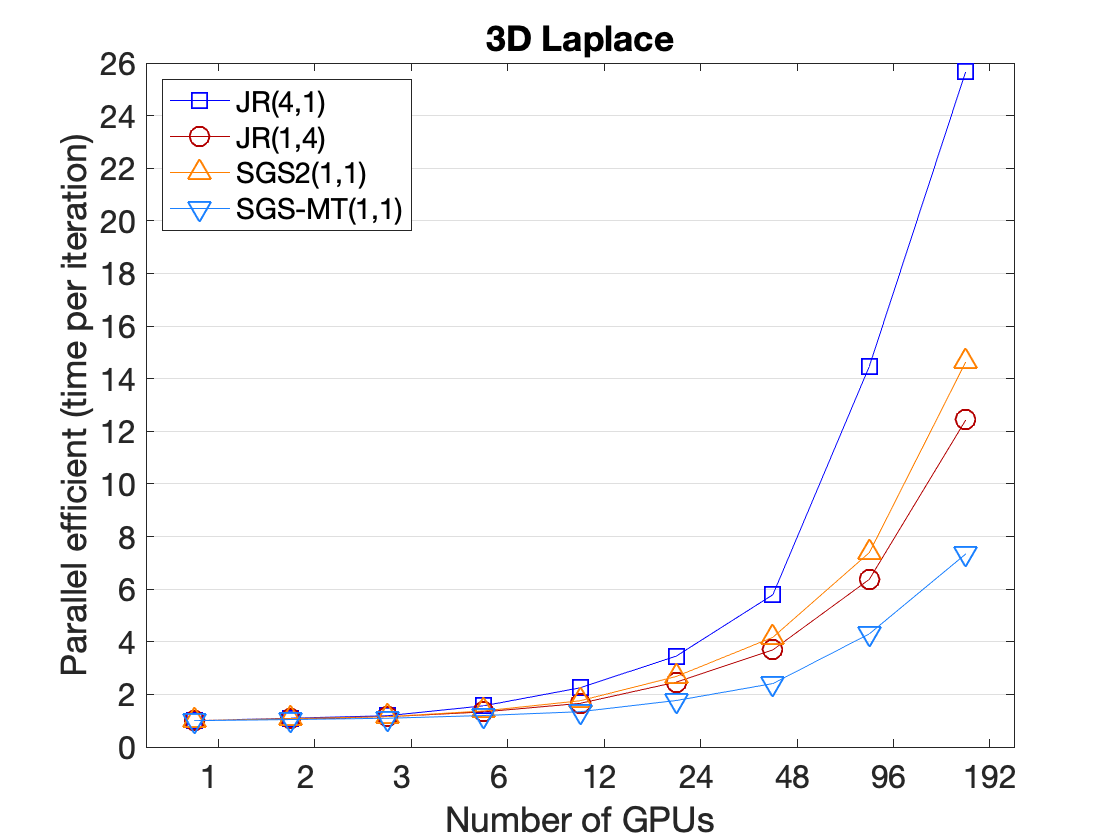}}
\caption{GMRES(60) parallel efficiency for Laplace 3D problems.}
\label{fig:laplace3d-parallel-efficiency}
\end{figure}

Finally, Figure~\ref{fig:laplace3d-strong-scale} displays the parallel 
strong-scaling when solving the Laplace 3D problem.
JR(4,1) and JR(1,4) are the ``global'' and ``local'' JR preconditioners
that apply the JR iteration on the global and local matrix, respectively.
The iteration count remains constant when using JR(4,1),
but may grow with the number of GPUs using JR(1,4) 
(see Figure~\ref{fig:laplace3d-parallel-efficiency}).
Thus, JR(4,1) performs well on a small number of GPUs 
(due to the faster solver convergence rate),
but JR(1,4) performs better on a larger number of GPUs 
(due to reduced communication cost).

%

%% file: result-ss.tex
Matrices from  the Suite-Sparse matrix collection are now used
to study the performance of the Trilinos SGS2 preconditioner.  
For these experiments,  symmetric positive definite (SPD) matrices 
are the focus and the conjugate gradient algorithm is applied
(CG)~\cite{Hestenes:1952} as our Krylov solver.  Figure~\ref{tab:suite-sparse}
lists the test matrices, and Figure~\ref{tab:cg-default} displays the performance.
When SGS2 requires just one inner sweep to obtain the same convergence rate as
that of the classical SGS, it can reduce the time to solution.  Also, SGS2
often outperformed JR(4), even though it may require more iterations
(because it was cheaper to apply).

For eight of the test matrices, with the default damping factor, 
JR diverged when using either the coefficient matrix or the triangular
part of the matrix.\footnote{In theory, JR on any triangular system 
should converge because the spectral radius of the iteration matrix, $\rho(D^{-1}L)$, is zero.
However, in our numerical experiments with SGS2,
inner JR failed to converge for some test matrices 
(e.g. a triangular matrix with a large condition number). 
Nevertheless, a carefully-selected damping factor
helped the convergence (e.g. reducing the condition number).}
Consequently, the solver did not converge using either JR(4) or
SGS2(1,1) for these matrices. 
%
%
Figure~\ref{tab:cg-damp} displays the
performance with SGS2(1,1) using a different damping factors~$\gamma$
for the inner JR sweeps.
The default damping factor was set for the outer sweep (i.e., $\omega=1.0$),
as in Figure~\ref{tab:cg-default}.
The robustness of the inner JR sweep,
and thus the convergence rate of SGS2 and GMRES, is improved
using a smaller damping $\gamma$.

\setlength{\tabcolsep}{3pt}
\begin{figure}
\tiny
\centerline{
 \begin{tabular}{l|rrrr|rrrrr|rrr|r}
         & \multicolumn{4}{c|}{CFD}
         & \multicolumn{5}{c|}{Structural Problem}
         & \multicolumn{3}{c|}{Thermal Problem}
         & \multicolumn{1}{c}{EM}\\
 \hline
id       & 1       & 2      & 3       & 4         & 5         & 6          & 7        & 8           & 9         & 10     & 11     & 12        & 13\\     
name     & PPoiss  &   cfd1 &    cfd2 & para\_fem & bcsstk17  & af\_shell7 & Fault639 & Emil923   & Hook1498  & ted\_B & therm1 & therm2    & tmt\_sym\\
 \hline
 $n$     &  14,822 & 70,656 & 123,440 & 525,825   & 10,974    & 504,855    & 638,802  & 923,136   & 1,498,023 & 10,605 & 82,654 & 1,228,045 & 726,713\\
 $\frac{nnz}{n}$ 
         &  48.3   & 25.8   & 25.0    & 7.0       & 39.1      & 34.8       & 42.7     & 43.7      & 39.6      & 13.6   & 6.9    & 7.0       & 7.0
 \end{tabular}
}
\caption{Test matrices (SPD) from Suite-Sparse matrix collection. 
Solution time (seconds). Rows $n$. Non-zeros per row $nnz(A)/n$.
\label{tab:suite-sparse}}
\end{figure}
\setlength{\tabcolsep}{6pt}

%
%
\begin{figure}
\scriptsize
\centerline{
 \begin{tabular}{l|r|rrrrr}
    & SGS                    & SGS-MT               & SGS2(1,1)           & JR(4)               &JR(1)\\
 \hline\hline
 1  &   1.59 ( 0.97,   243)  &  0.78 ( 0.70,   269) & $--$                & $--$                & $--$\\
 2  &   8.30 ( 4.99,   486)  &  1.57 ( 1.33,   648) & $--$                & $--$                & $--$\\
 3  &  54.25 (33.29, 1,945)  &  8.04 ( 6.96, 2,577) & $--$                & $--$                & $--$\\
 4  &  57.69 (37.35, 1,434)  &  6.39 ( 5.34, 1,465) & 2.71 ( 1.69, 1,448) & 4.18 ( 3.43, 1,050) & 3.25 (1.78, 2,099)\\
 \hline
 5  &   4.93 (  2.90, 1,122) &  3.64 ( 3.25, 1,160) & $--$                & $--$                & $--$\\
 6  & 212.20 (130.70, 1,400) & 17.83 (16.61, 1,407) & $--$                & $--$                & $--$\\
 7  & 831.30 (528.40, 3,339) & 56.48 (53.76, 2,839) & $--$                & $--$                & $--$\\
 8  & 1089.0 (683.30, 3,029) & 100.6 (96.05, 3,680) & $--$                & $--$                & $--$\\
 9  &  989.5 (652.40, 1,785) & 81.66 (78.13, 1,972) & $--$                & $--$                & $--$\\
 \hline
 10 & 0.033  (0.015, 11)     & 0.026 (0.021, 14)    & 0.015 (0.006, 16)   & 0.329 (0.249, 329)   & 0.269 (0.147, 505) \\
 11 & 5.54   (3.2, 701)      &  0.87 (0.61, 726)    & 0.55  (0.34, 718)   & 0.541 (0.385, 408)   & 0.524 (0.213, 821) \\
 12 & 278.50 (196.50, 2,657) & 24.54 (21.81, 2,731) & 7.67  (4.98, 2,714) & 9.27  (7.75, 1,541)  & 7.05  (4.01, 3,090)\\
\hline
 13 & 105.60 (69.67, 1,857)  & 13.30 (11.44, 2,363) & 4.41  (2.74, 2,139) & 23.22 (19.18, 5,235) & 16.63 (9.18, 9,632)
 \end{tabular}
}
\caption{\label{tab:cg-default}
Solution time in seconds with Suite-Sparse matrices. The numbers in 
parentheses are time to apply the preconditioners and CG iteration counts.}
\end{figure}

\setlength{\tabcolsep}{3pt}
\begin{figure}
\tiny
\centerline{
\begin{tabular}{l|ccccccccc}
         & \multicolumn{9}{c}{Matrix ID}\\
$\gamma$ & 1             & 2              & 3              & 5             & 6            & 7             & 8             & 9             & 12\\
\hline 
0.1      & 0.47 (626)    &  1.47 (1,627)  &  7.09 (6,606)  & 1.71 (2,342)  & 7.59 (3,147) & 21.13 (6,485) & 39.09 (8,833) & 30.64 (4,708) & 10.24 (4,822)\\
0.2      & 0.38 (503)    &  1.11 (1,221)  &  5.13 (4,823)  & 1.30 (1,782)  & 6.91 (2,845) & 17.52 (5,376) & 32.50 (7,351) & 26.42 (4,069) &  9.27 (4,347)\\
0.3      & 0.31 (415)    &  0.91 (998)    & 16.51 (15,535) & 1.07 (1,474)  & 5.89 (2,429) & 14.83 (4,555) & 27.31 (6,181) & 22.93 (3,529) &  8.38 (3,921)\\
0.4      & 0.27 (359)    & 10.49 (11,633) & $--$    ($--$) & 3.39 (4,657)  & 5.14 (2,120) & 12.65 (3,878) & 23.64 (5,346) & 20.07 (3,085) &  7.58 (3,544)\\
0.5      & 0.33 (456)    & $--$  ($--$)   & $--$    ($--$) & $--$  ($--$)  & 4.59 (1,894) & 11.38 (3,491) & 21.12 (4,769) & 17.93 (2,755) &  6.84 (3,210)\\
\end{tabular}}
\caption{CG iteration counts with SGS2(1,1) preconditioner. 
The number in parentheses displays the CG iteration count.}
\label{tab:cg-damp}
\end{figure}
\setlength{\tabcolsep}{6pt}


%% file: result-mixed.tex
\begin{figure}[h]
\centerline{
  \includegraphics[width=.5\textwidth]{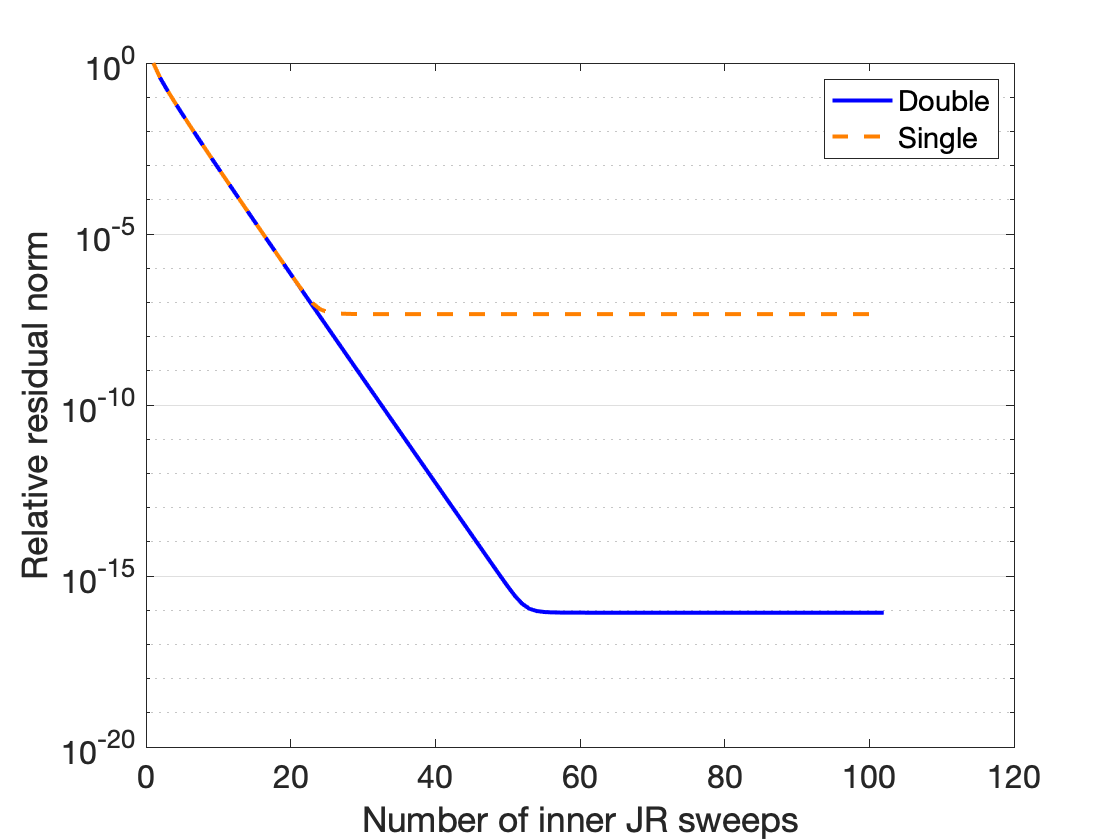}
  \includegraphics[width=.5\textwidth]{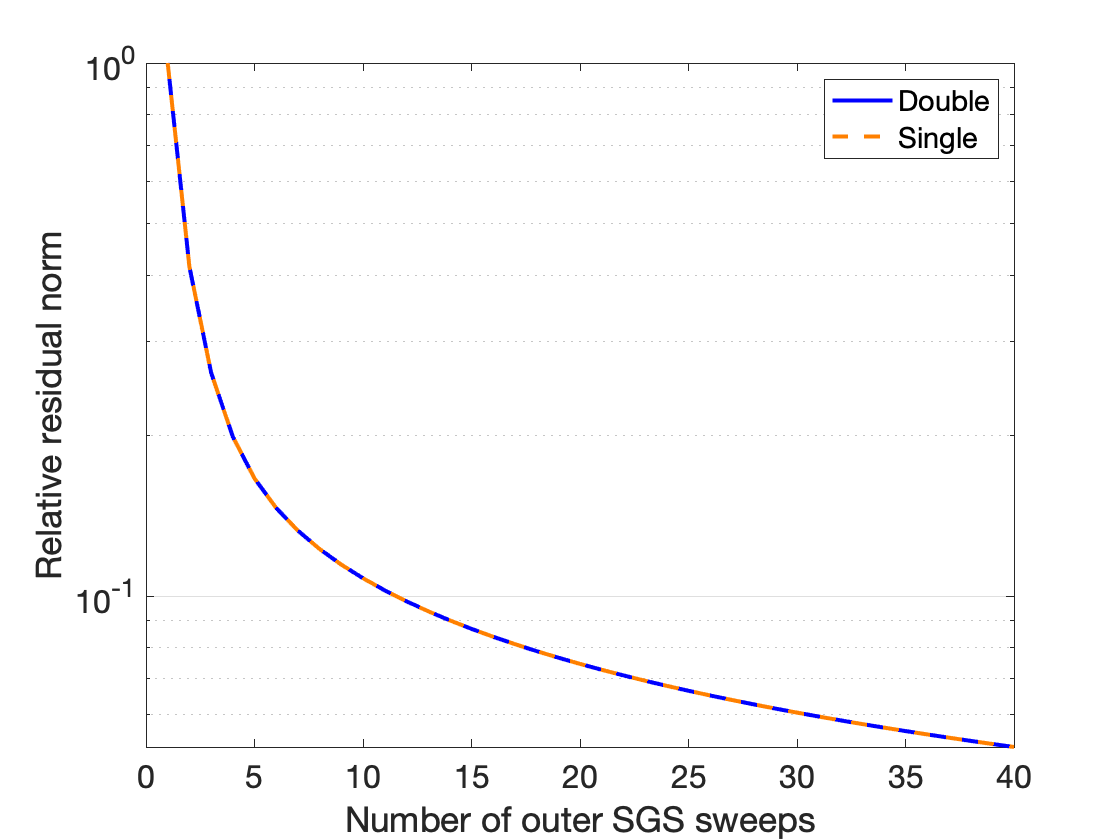}}
\caption{\label{fig:lower-prec}
Outer SGS2(1,1) and Inner JR convergence for 2D Laplace ($n_x = 1000$).}
\end{figure}

\begin{figure}[h]\scriptsize
\centerline{
\subfloat[GMRES(60).]{
\begin{tabular}{l||rr|rr|rr}
                  & SGS, double & SGS, single & SGS2, double & SGS2, single & MT, double, & MT, single\\
\hline\hline
\# of GMRES iters & 5,095       & 5,095       & 6,656        & 6,656        & 8,236       & 7,745\\
\hline
time to apply SGS & 201.9       & 178.1       & 8.13         &  5.63        & 28.28       & 25.99\\
time to solution  & 315.7       & 205.5       & 23.21        & 21.33        & 46.68       & 43.99\\
\end{tabular}}}
\centerline{
\subfloat[CG.]{
\begin{tabular}{l||rr|rr|rr}
                  & SGS, double & SGS, single & SGS2, double & SGS2, single & MT, double, & MT, single\\
\hline\hline
\# of CG iters    & 1,108       & 1,108       & 1,279        & 1,279        & 1,627       & 1,911\\
\hline
time to apply SGS & 43.07       & 38.01       & 1.54         &  1.07        & 5.48        & 6.28\\
time to solution  & 65.35       & 41.94       & 2.49         &  2.13        & 6.71        & 7.77\\
\end{tabular}}}
%
\caption{Compute time (seconds)
of CG with 2D Laplace $(n_x =1000)$ using SGS preconditioner in
double or single precision.\label{fig:lower-prec-laplace}}
\end{figure}

In experiments with the GS2 preconditioner,
only one outer and inner sweep were enough to obtain a convergence
rate similar to that with the classical GS preconditioner.
Hence, high-accuracy is not needed from the outer GS iterations or 
from the inner JR iterations.
To improve the performance,  the use of
lower-precision floating point arithmetic was investigated 
for the preconditioner.
Our main concern is to determine if the quality of the
preconditioner would degrade in lower precision.  This
depends on the properties of the matrix that affect the behavior of GS in
finite precision, such as its maximum attainable accuracy and convergence rate.
For example, GS is norm-wise backward stable and its attainable relative residual norm
is upper-bounded by $\epsilon \kappa(A)$~\cite{Higham::1993}, while each GS or JR sweep would reduce
the relative residual norm by a factor of at least $\|I - M^{-1}A\|$.
Using single precision for GS, the quality of the preconditioner
should be maintained (same as GS in double precision)
if one sweep of GS in double precision does not
reduce the residual norm below the attainable accuracy of GS in single precision
(e.g. $\kappa(A) < \mathcal{O}(1/\epsilon_s)$ and 
$\mathcal{O}(\epsilon_s) \kappa(A) < \|I - M^{-1}A\|$,
where $\epsilon_s$ is the machine precision in single precision).

Figure~\ref{fig:lower-prec} indicates that the outer GS and inner JR iteration
in single precision converge in the same way as in double precision until the iteration achieves the maximum solution accuracy and begins to stagnate.
Figures~\ref{fig:lower-prec-laplace} and \ref{fig:lower-prec-ss} 
then display the performance of CG using the SGS preconditioner 
in double or single precision. The CG convergence rate was maintained 
using the single-precision preconditioner for these test matrices, 
while the time for applying the preconditioner was reduced by factors
of around $1.4\times$.

\setlength{\tabcolsep}{3pt}
\begin{figure}[h]\tiny
\begin{tabular}{l||rrr|rrr|rrr|rrr}
                  &\multicolumn{3}{c|}{4  ($\gamma=1.0$)}
                  &\multicolumn{3}{c|}{6  ($\gamma=0.5$)}
                  &\multicolumn{3}{c|}{7  ($\gamma=0.5$)}
                  &\multicolumn{3}{c}{12 ($\gamma=1.0$)}\\
                  & double & single & speedup & double & single & speedup & double & single & speedup\\
\hline\hline
\# of CG iters    & 1,448  & 1,448  & --      &  1,894 & 1,894  & --      & 3,491  & 3,532  & --      & 2,714  & 2,714  & --\\
\hline
time to apply SGS &  1.69  & 1.14   & 1.48    & 3.30   & 2.33   & 1.41    &  8.01  & 5.99   & 1.35    & 4.98   & 3.55   & 1.40\\
time to solution  &  2.71  & 2.34   & 1.16    & 4.59   & 3.86   & 1.18    & 11.38  & 9.75   & 1.17    & 7.67   & 6.89   & 1.11\\
\end{tabular}
\caption{Compute time (seconds)  CG with the Suite Sparse matrices in Figure~\ref{tab:suite-sparse}
using SGS2 preconditioner in double or single precision.\label{fig:lower-prec-ss}}
\end{figure}
\setlength{\tabcolsep}{6pt}

%% file: result.tex
To study the performance of the two-stage Gauss-Seidel preconditioner
and smoother in a practical setting, incompressible fluid flow 
simulations were performed with Nalu-Wind~\cite{Sprague:2020}.
This is the primary fluid mechanics code for the ExaWind project,
one of the application projects chosen for the DOE Exascale Computing Project (ECP)
and is used for high-fidelity simulations of air flow dynamics around wind turbines.
Nalu-Wind solves the acoustically incompressible Navier-Stokes
equations, where mass continuity is maintained by an approximate pressure
projection scheme.  The governing physics equations for momentum, pressure, and
scalar quantities are discretized in time with a second-order BDF-2 integrator,
where an outer Picard fixed-point iteration is employed to reduce the nonlinear
system residual at each time step.
Within each time step, the Nalu-Wind simulation time is often dominated by the
time required to setup and solve the linearized governing equations,
using either \hypre or Trilinos solver package.
To solve the momentum equations, Nalu-Wind typically employs Gauss-Seidel (GS) or symmetric
Gauss-Seidel (SGS) iteration as a preconditioner to accelerate the GMRES convergence.
The pressure systems are solved using GMRES with an
algebraic multigrid (AMG) preconditioner, where a Gauss-Seidel or Chebyshev
polynomial smoother is applied to relax or remove high energy components of the
solution error (e.g. those associated with the large eigenvalues of the
system), which the coarse-grid solver fails to eliminate.  
%

To study the performance of GS2, representative
incompressible Navier-Stokes wind-turbine simulations were performed 
(an atmospheric boundary layer (ABL) precursor simulation on a structured mesh).
%
The strong scaling performance of
both the Trilinos and \hypre solver stacks were evaluated.
We note that it is possible to generate segregated linear equations for solving
the momentum problems. This reduces the sizes of the linear systems,
and also reduces the number of nonzero entries per row.
This technique can reduce the time to solutions, especially when the underlying
computation kernels are optimized to take advantage of the multiple right-hand-sides.
In the following, this option was not used for the experiments using Trilinos, while this option was employed for the experiments using \hypre.



\subsection{Experimental Results with Trilinos}

\input{result-smooth}
\input{result-trilinos}

\subsection{Experimental Results with \hypre}
\input{result-hypre}

%% file: result-trilinos.tex
\begin{figure}
\centerline{
 \subfloat[For an SGS sweep of the momentum solves.]{
  \includegraphics[width=.5\textwidth]{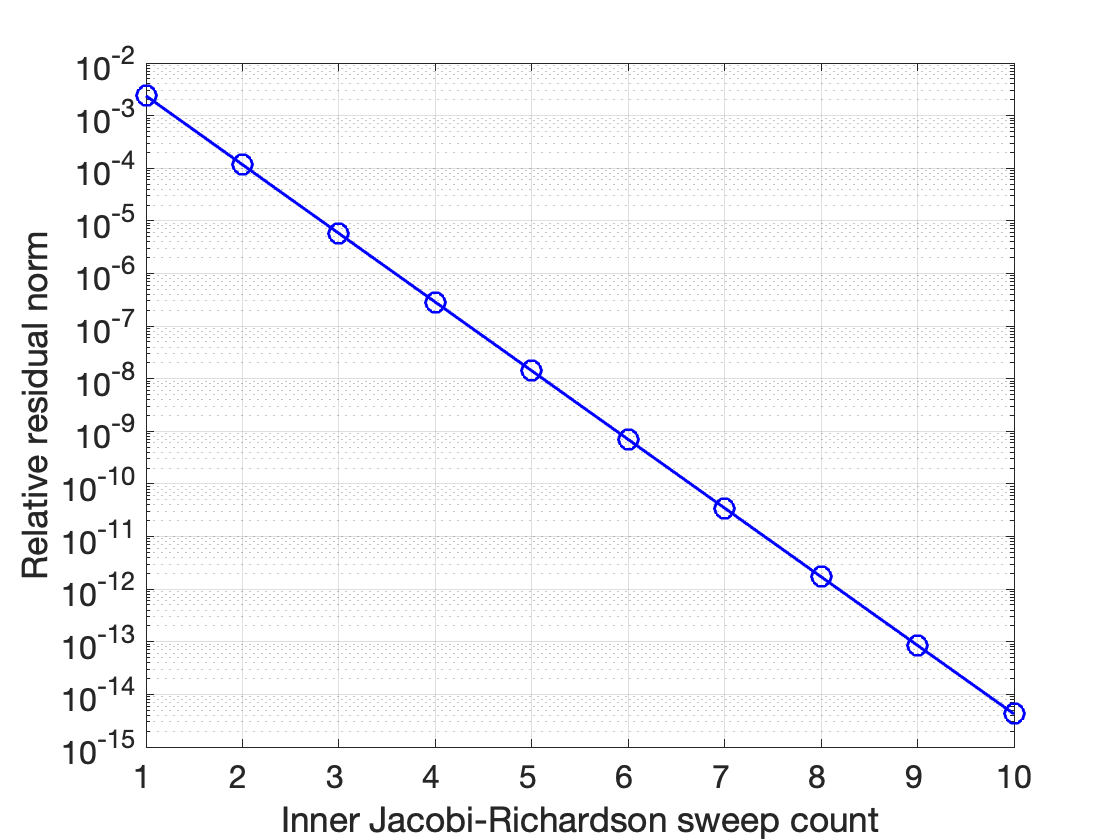}}
 \quad
 \subfloat[For an SGS sweep of the  continuity solve.]{
  \includegraphics[width=.5\textwidth]{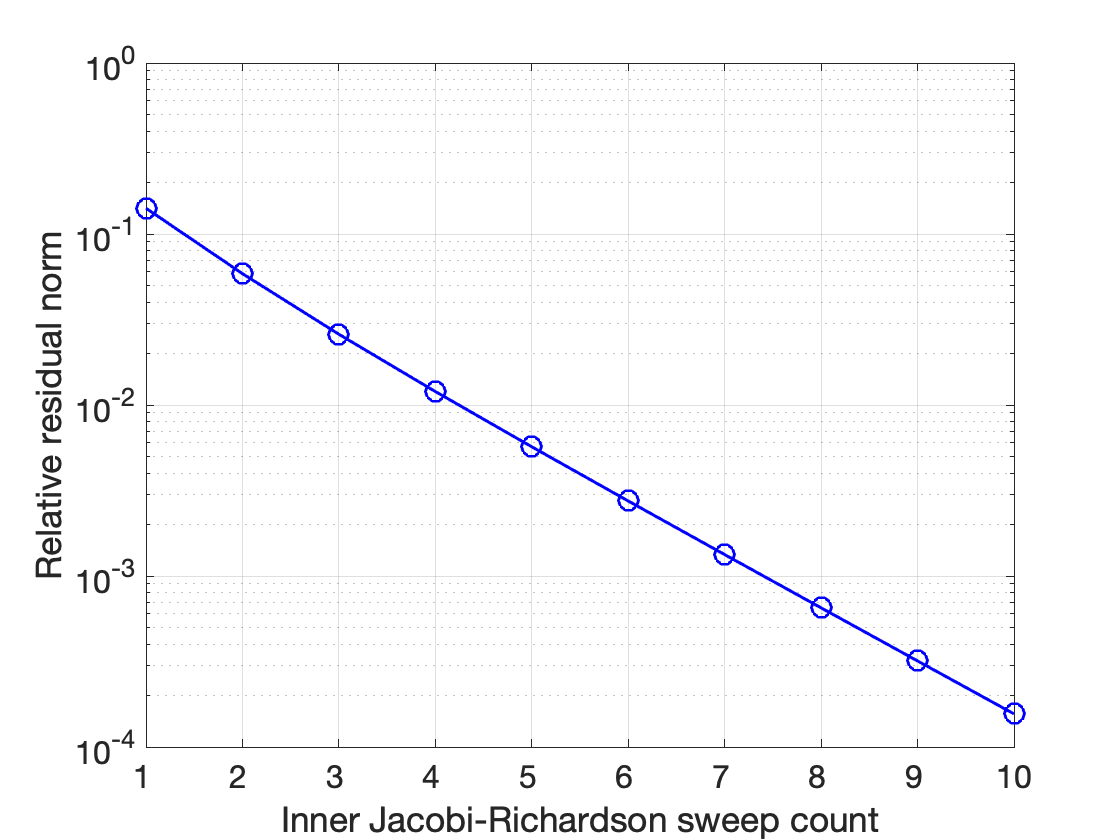}}
}
\caption{\label{fig:inner-convergence}
Convergence history of the relative residual norm, $\|\vect{r}_l -
(L+D)\vect{g}_k^{(j)}\|_2/\|\vect{r}_k\|_2$, with inner Jacobi-Richardson
sweeps for the ABL 40m simulation.
}
\end{figure}

Figure~\ref{fig:inner-convergence} displays
the convergence history of the inner sweep for one time-step of the Nalu-Wind
ABL precursor simulation on a 40m resolution mesh.
In this simulation, the solver converged with the same number of iterations
to a relative residual tolerance level of $1\times 10^{-6}$
using just one sweep of the inner Jacobi iteration for SGS2
as that using the sequential SGS preconditioner.
%
Moreover, Figure~\ref{tab:compact} reports the iteration counts for ten time-steps 
of the Nalu-Wind ABL precursor 40m simulation,
where the compact form leads to a significantly larger number of iterations
for the pressure systems.
On another note, the compact form worked well as the preconditioner for the remaining problems
(with an extra inner sweep for the compact form, the solver often converged similar to 
with the non-compact form).
A comparison of these is presented in Section~\ref{sec:others}
along with similar results for model problems in Section~\ref{sec:result-model}.

\begin{figure}
\begin{center}
\begin{tabular}{l||r|r|r||rr|rrrr}
            &     &      &  & \multicolumn{5}{c}{two-stage}\\
            & No  & Seq  & MT  & \multicolumn{2}{c|}{non-compact} & \multicolumn{4}{c}{compact}\\
\hline
Inner sweep & --  & --  & --   & 0   & 1   & 1    & 2    & 4    & 6\\
Continuity  & --  & 285 & 311  & 511 & 300 & 2410 & 863  & 401 & 314\\
Enthalpy    & 419 & 81  & 121  & 122 & 84  & 84   & 82   & 81  & 81\\
Momentum    & 270 & 83  & 108  & 81  & 83  & 90   & 83   & 83  & 83\\
TKE         & 334 & 80  & 120  & 104 & 81  & 81   & 80   & 80  & 80\\
\end{tabular}
\caption{Solver iteration counts with ten time-step 
ABL 40m simulation on one GPU with SGS preconditioner 
for momentum and turbulent kinetic energy (TKE) systems 
(one SGS sweep) and continuity (three SGS sweeps).
\label{tab:compact}
}
\end{center}
\end{figure}

Next, consider the performance of Trilinos solvers for the Nalu-Wind simulation
using the ABL precursor 20m mesh.  To solve the continuity system, Trilinos
uses a smoothed aggregation SA-AMG preconditioner and a Chebyshev
smoother.  The three remaining equations (Momentum, TKE and Enthalpy) are
solved with a GS preconditioner.  
In this section, GS2
is applied as alternative to the sequential GS preconditioner.
%
\ignore{
Figure~\ref{tab:smooth-trilinos} displays the number of iterations
for solving the continuity problems, by varying the number of GS smoother sweeps
and damping factors. For this purpose, three different combinations of 
GS smoothers were applied: 
1) forward GS for both pre and post smoothers,
2) forward GS for the pre smoother and backward GS for the post smoother, and
3) SGS for both pre and post smoothers, and the second option often
leads to the lowest iteration counts.
For this problem, especially with the carefully-selected damping factors, the
Chebyshev and JR smoothers performed well, although similar trends are observed
when just one inner sweep is employed. When compared with the GS2 smoother,
the solver converged with slightly more iterations using JR(3) but needed fewer
iterations using JR(4).
}

\begin{figure}[!t]
\centerline{
\subfloat[Strong scaling Trilinos.] {\label{trilinos-strong}\includegraphics[width=.6\textwidth]{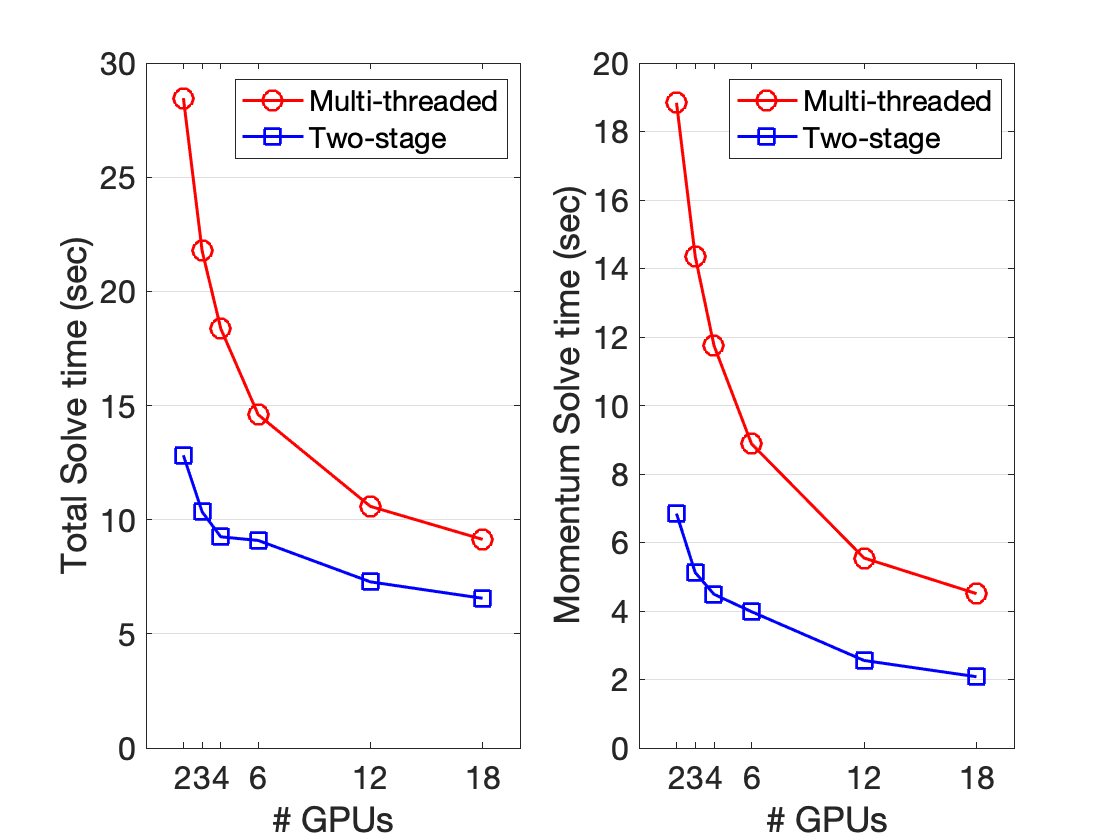} }
\subfloat[Strong scaling \hypre.] {\label{hypre-strong}\includegraphics[width=.6\textwidth]{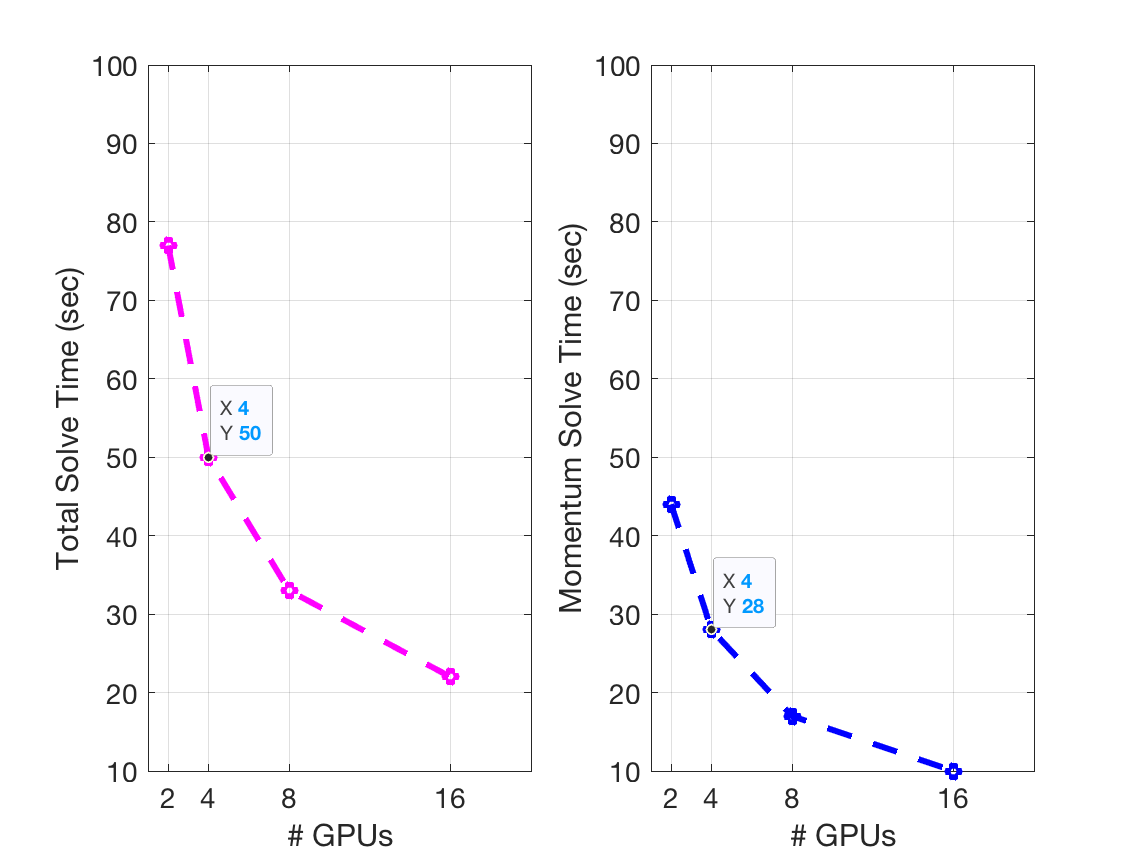} }
}
\caption{Parallel scalability of solve time on Summit for ABL 20m precursor.} 
\end{figure}

Figure \ref{trilinos-strong} 
displays the
strong scaling of the Trilinos solvers for the Nalu-Wind
simulation using the ABL precursor mesh.
The maximum number of iterations required to solve these linear systems was low 
but remains constant across all mesh resolutions tested.
Even though the same Chebyshev smoother was used for solving the pressure system, 
compared to using the classical SGS preconditioner, the total solve time was 
still reduced using SGS2 preconditioner for the remaining systems. The performance
advantage of using the SGS2 preconditioner over the multi-threaded SGS was 
maintained over multiple nodes.


%% file: result-hypre.tex
Strong scaling studies for the \hypre solver stack have also been
performed for the ABL precursor Nalu-Wind simulation.  
The GS2 iteration was implemented as both the
momentum equation preconditioner and as a smoother for the GMRES solver with
classical Ruge-St\"uben C-AMG preconditioner for the pressure continuity solver
\cite{Ruge87}. 
For the C-AMG set-up, HMIS coarsening
is specified for the ABL precursor problem.
One sweep of the GS2 smoother is applied at
each level of the C-AMG $V$-cycle. The strength of connection threshold is
$\theta = 0.25$ and two levels of aggressive coarsening reduce the
preconditioner cost and coarse-grid operator complexity.

MGS-GMRES solver compute times are 
reported for ten time-step integrations from a restart of the Nalu-Wind model
~\cite{Thomas:2019, Swirydowicz:2020} .
Momentum solve time exhibits nearly perfect strong scaling,
every time the number of ranks is doubled, the compute time is halved (Figure
\ref{hypre-strong}, right panel). The continuity solver time
exhibits sub-linear strong scaling. This may also be attributed to additional
computation and communication overheads in the C-AMG $V$-cycle such as the MPI
communication associated with the prolongation and restriction matrix-vector
products. The authors are currently investigating the source of these
sub-linear scaling solver times.

%% file: conclusion.tex
A two-stage iterative variant of the Gauss-Seidel relaxation
was examined as an alternative to the sequential GS 
algorithm based on a sparse triangular solver. The two-stage approach replaces
the triangular solve with a fixed number of inner relaxation sweeps, and
is much more amenable to an efficient parallel implementation on many-core architectures. 
When a small number of inner sweeps is needed,
the two-stage variant often outperforms the sequential approach.


To study the performance of the GS2 algorithm  as a preconditioner and smoother 
in a practical setting, the algorithm was employed 
in Nalu-Wind incompressible fluid flow simulations~\cite{Sprague:2020}.
GS2 was applied to accelerate convergence for the momentum equation and as a 
smoother for the AMG preconditioner when for the pressure solver.  Our
performance results demonstrate that the GS2 preconditioner
requires just one inner sweep to obtain a GMRES convergence rate
similar to the sequential algorithm, and thus it achieves a faster time to
solution, both on CPU and GPU.
Another advantage, not discussed earlier, is that the set-up cost for 
GS2 is minimal and does not require eigenvalue estimates.  
The setup cost influences the choice of preconditioner, because the setup
is performed at each time step due to mesh movement.
In addition, to obtain the desired convergence rate of SGS2 and GMRES,
the required parameters (e.g., number of inner sweeps 
and damping factors) depend on the properties of the coefficient matrices.
The relationship between the inner stopping criteria
and the outer convergence rate is being examined in order to determine
the appropriate stopping criteria.